\documentclass[reqno,final]{amsart}
\usepackage{natbib}  
\usepackage{fancyhdr} 
\usepackage{color} 
\usepackage{hyperref} 
\usepackage{graphicx} 

\usepackage{pstricks}
\usepackage{amssymb}
\usepackage{lineno}

\definecolor{aleacolor}{rgb}{0.16,0.59,0.78}

\hypersetup{
breaklinks,
colorlinks=true,
linkcolor=aleacolor,
urlcolor=aleacolor,
citecolor=aleacolor}

\renewcommand{\cite}{\citet}

\theoremstyle{plain}
\newtheorem{theorem}{Theorem}[section]                                          
\newtheorem{proposition}[theorem]{Proposition}                          

\newtheorem{corollary}[theorem]{Corollary}

\theoremstyle{definition}

\theoremstyle{remark}

\makeatletter \@addtoreset{equation}{section} \makeatother

\newcommand{\Di}{\displaystyle}

\begin{document}

\title[Expansions for Sample Median with Random Sample Size]{Second Order Expansions for Sample Median with Random Sample Size}

\author{Gerd Christoph}
\author{Vladimir V. Ulyanov}
\author{Vladimir E. Bening}

\address{Otto-von-Guericke University Magdeburg, Department of Mathematics, \newline
Postfach 4120,\newline
39016 Magdeburg, Germany.}
\email{gerd.christoph@ovgu.de}

\address{Lomonosov Moscow State University,\newline Faculty of Computational Mathematics and Cybernetics\newline 
119991, Leninskie Gory, 1/52, Moscow, Russia.\newline $\quad$
National Research University Higher School of Economics,\newline
101000,  Myasnitskaya ulitsa, 20, 
 Moscow, Russia}
\email{
vulyanov@cs.msu.ru}

\address{Lomonosov Moscow State University,\newline Faculty of Computational Mathematics and Cybernetics\newline 
119991, Leninskie Gory, 1/52, Moscow,  Russia}
\email{bening@yandex.ru}


\subjclass[2000]{60F05, 60G50, 62E17, 62H10.} 
\keywords{Sample median; samples with random sizes; second order expansions; Laplace distribution; 
Student's $t$-distribution; negative binomial distribution; discrete Pareto distribution.}

\begin{abstract}
 In practice, we often encounter situations, where a sample size is not defined in advance and can be a random value. 
The randomness of the sample size crucially changes the asymptotic properties of the underlying statistic. 
In the present paper the second order Chebyshev--Edgeworth and Cornish--Fisher expansions based of Student's $t$- and 
Laplace distributions and their quantiles are derived for sample median with random sample size of a special kind. 
\end{abstract}

\maketitle

\section{Introduction}
Usually in classical statistical inference the number of observations is known. But often we do not know in advance the sample sizes or there 
are missing observations. Therefore the sample size may be a realization of a random variable.

There are many practical situations, where it is almost impossible to have a fixed sample size. They often occur
when observations are collected in a fixed time span. For example, in reliability testing this is the number of 
failed devices, in medicine -- the number of patients with a specific disease, in finance -- the number of market transactions, in queueing theory  -- the number of customers entering a store,
in insurance -- the number of claims. All these numbers are random variables.

The use of samples with random sample sizes has been steadily growing over the years. 
For an overview of statistical inferences with a random number of observations and some applications see, e.g. \citet{EMM16} and the references therein.

Let $X_1, X_2, \ldots \in \mathbb{R}=(- \infty\,,\, \infty)$  and $N_{1}, N_{2}, \ldots \in \mathbb{N}=\{1, 2, ...\}$ be the  random variables on the 
same probability space $\left(\Omega,\mathbb{A},\mathbb{P}\right)$. In statistics the random variables 
$X_1, X_2, \ldots $ are observations.
 Let $N_n$ be a  
random size of the underlying sample, which depends on  parameter $n \in \mathbb{N}$. We suppose for each $n \in \mathbb{N}$ that  $N_n \in  \mathbb{N}$   
is independent  of  $X_1, X_2, \ldots$ and   $N_n \to \infty$  in probability as $n \to \infty$.

Let   $T_m := T_m \left(X_1, \ldots, X_m\right)$ be some statistic of a sample with non-random sample size $m \in \mathbb{N}$.
Define the random variable $\overline{T}_{N_n}$ for every $n \in \mathbb{N}$:
\begin{equation}\label{eq0} 
\overline{T}_{N_n} (\omega) := \overline{T}_{N_n(\omega)} \left(X_1 (\omega), \ldots, X_{N_n (\omega)} \right),\,\,\,\omega \in \Omega,
\end{equation}
i.e. $\overline{T}_{N_n}$ is some statistic obtained from a random  sample $X_1, X_2, \ldots, X_{N_n}$.

\citet{Gn89} considered the asymptotic properties
of the distributions of sample quantiles for  samples of random size.
In \citet{Nun19a} unknown sample sizes are assumed  in medical research for analysis of one-way fixed effects ANOVA 
models to avoid false rejections.  Application of orthogonal mixed models to situations with sample of random sizes 
are investigated in \citet{Nun19b}.
\citet{EMM16} considered  inference for the mean with known and unknown
variance and inference for the variance in the normal model.
Prediction intervals for the future observations for generalized order statistics and confidence intervals for quantiles based on 
samples of random sizes are studied in \citet{BNEY18} and \citet{AMR17}, respectively. They illustrated  
their results with real biometric data set, the duration of remission of leukemia patients treated by one drug.
General asymptotic expansions for statistics with random sample sizes $\overline{T}_{N_n}$ are given in \citet{BGK13} applying corresponding asymptotic expansions 
for the normalized statistic $T_m$ and the suitable scaled random sample size $N_n$.

Many models lead to random sums and random means
\begin{equation}\label{eq01}
S_{N_n} = \sum\nolimits_{k=1}^{N_n} X_k \quad \mbox{and} \quad T_{N_n}  = \frac{1}{N_n}\sum\nolimits_{k=1}^{N_n} X_k= \frac{1}{N_n}S_{N_n},
\end{equation}
respectively. 
Wald's  identity for random sums $\mathbb{E}(S_{N_n})=\mathbb{E}(N_n) \mathbb{E}(X_1)$ if $N_n$ and $X_1$ 
have finite expectations is a powerful tool in statistical inference, particularly in sequential analysis, 
see e.g. \citet{Wa45} and  \citet{KP49}.
\citet{Rob48} proved that asymptotic normality of the index $N_n$ automatically implies asymptotic normality 
of the corresponding random sum $S_{N_n}$.

The randomness of  the sample size may crucially change asymptotic properties of random sums, 
see e.g. \citet{Gn89} or \citet{GK96}. If the statistic $T_m$ is asymptotically normal, then
the limit laws of normalized statistic $T_{N_n}$ are scale mixtures of normal distributions with zero mean, depending on the 
random sample size $N_n$.

A fundamental introduction to asymptotic distributions of random sums is given in \citet{Doe15}.
Using Stein's method,  
quantitative Berry-Esseen bounds  of random sums  were proved in \citet[Theorem 10.6]{CGS11}, 
\citet[Theorems 2.5 and 2.7]{Doe15} and \citet[Theorem 1.3]{PR14} in case of approximation  by  normal and Laplace distributions.
Moderate and large derivations are investigated in \citet{EL19}, and \citet{KM97}. 
Many applications of geometric random sums when $N_n$ is geometrically distributed are given in \citet{Ka97}.
Bounds on the total variation distance between geometric random sum of independent, non-negative,
integer-valued random variables and the geometric distribution are studied in \citet[Section 3]{PRR13}

It is worth to mention that a suitable scaled factor by random sums $S_{N_n}$ or random means $T_{N_n}$ affects the type of limit distribution.
In fact, consider random mean $T_{N_n}$  given in (\ref{eq01}). For the sake of convenience let $X_1,X_2,... $ be 
independent standard normal random variables and $N_n \in \mathbb{N}$ be geometrically distributed with 
$\mathbb{E}(N_n)=n$ and independent of $X_1, X_2,...$. Then one has 
\begin{eqnarray}
\bullet \,&  \mathbb{P}\left(\sqrt{\Di N_n}\,T_{N_n} \leq x\right)
 \quad = \quad \int\limits_{-\infty}^x \frac{\Di 1}{\Di \sqrt{2\,\pi}} e^{- u^2/2} du \quad \mbox{for all} \,\, n \in \mathbb{N},\label{sum1}\\
\bullet  \,&  \mathbb{P}\left(\sqrt{\mathbb{E}(N_n)}\, T_{N_n} \leq x\right)
 \quad \to \quad \int\limits_{-\infty}^x  \left( 2 +  u^2\right)^{-3/2} du  \quad \mbox{as} \,\, n \to \infty . \label{sum3}\\
\bullet \,&  \mathbb{P}\left(\frac{\Di N_n}{\Di \sqrt{\mathbb{E}(N_n)}}\, T_{N_n}\leq x\right)
  \quad \to \quad \int\limits_{-\infty}^x  \frac{\Di 1}{\Di \sqrt{2}} e^{- \sqrt{2}\,|u|} du  \quad \mbox{as} \,\, n \to \infty,\label{sum2}
\end{eqnarray} 

\noindent
We have three different limit distributions. The suitable scaled random mean $T_{N_n}$  is 
standard normal distributed or tends to the Student distribution with 2 degrees of freedom
as the limit distributions depending on whether we take the random scaling factor $\sqrt{N_n}$ or  
the non-random scaling factor $\sqrt{\mathbb{E}N_n}$, respectively. Moreover, we get the Laplace distribution with variance 1  if we use scaling with the mixed  factor $N_n/\sqrt{\mathbb{E}(N_n)}$.

Assertion (\ref{sum1}) we obtain by conditioning and the stability of the normal law. Student distribution as a limit for  statistics from samples with a random sample size are proved e.g. in \citet{BK05} and \citet{ST16}, hence relationship 
(\ref{sum3}) holds. Since $N_n \, T_{N_n} = S_{N_n}$ statement (\ref{sum2})
follows e.g. from \citet{BK08} or \citet{ST16}.

In \citet{BGK13} first order expansions of the random mean $T_{N_n} = (X_1+...+X_{N_n})/N_n$ are proved if the sample 
size $N_n$ is negative binomial distributed with success probability $1/n$ or it is the maximum of $n$ independent 
identically distributed discrete Pareto random variables with tail index 1, using first order Chebyshev-Edgeworth 
expansions for mean $T_m = (X_1+...+X_m)/m$ and the rate of convergence for the distribution of suitably normalized random sample size $N_n$ to the corresponding limit law. Second order asymptotic expansions  of suitably normalized random sample size $N_n$ are proved in \citet{CMU20} which were used to derive second order Chebyshev-Edgeworth expansions for the random mean $T_{N_n} = (X_1+...+X_{N_n})/N_n$.\\[1ex]

In the present paper we investigate the median of a sample $\{X_1,....,X_{N_n}\}$ with the  
random sizes $N_n$ mentioned above.

Let  $F_X(x \, - \, \theta)$ and  $p_X(x \, - \, \theta)$ be the known common distribution function and the probability 
density function of independent components of the sample $\{X_1, X_2, ..., X_m\}$, where $\theta$ is the unknown location
parameter to be estimated from the given sample.
By $X_{(1)} \leq X_{(2)} \leq ... \leq X_{(m)}$ we denote the order statistics
constructed from the original observations $X_1, X_2, ..., X_m$.

As statistic $T_m$ we consider the sample median $M_m$, that is,
\begin{equation}\label{e1}
 M_m= \left\{ \begin{array}{ll} X_{(j)}, & \quad m=2 j - 1 , \\[1ex]
(X_{(j)}+X_{(j+1)})/2, & \quad m=2 j , \end{array} \right. \qquad j, m \in \mathbb{N}\,.
\end{equation}

\citet{Hu99} discussed the even-odd phenomenon for the median in statistical literature and gave a counterexample which contradicts the statistical
folklore: ``It never pays to base the median on an odd number of observations''.

Looking for change points in the location parameter in time series, tests for a 
change in mean may be susceptible to outliers in the data, whereas tests for a change
in median could may show a change of the center of the marginal distribution,  see \citet{SZ10}, 
\citet{VW17} and the references therein.

To perform statistical analysis of large data sets \citet{Mi19} presents new results for the {\it median-of-means estimator} 
using new algorithms for distributed  statistical estimation that exploit {\it divide-and-conquer approach}.

To estimate the location parameter $\theta$ one could use the random mean $T_{N_n}$ as well, but for its  second order  expansion more than the fourth moment of $X_1$ is required.   For heavy tailed distributions $F$ of $X_1$ with tail index $ \leq 4$ such second order Edgeworth expansions of the random mean $T_m$ cannot be obtained. If the tail index $\leq 1$, then the  mean does not exist: $\mathbb{E}|X_1|= \infty$. The mean need not always exist, whereas the median always exists. 

In \citet{PK19} confidence region for median of $X_1$ in  the nonparametric measurement error model are constructed and several  applications are given when a confidence interval about the center of a distribution is desired.

Therefore, it is reasonable to use the sample median $M_m$. 

The asymptotic normality of the normalized sample median
$M_m$ is well known, see e.g. \citet[Chapter 28.5]{Cr46}: 
If $F_X(0) = 1/2$, $p_X(0) > 0$ and  the density $p_X(x)$ is continuous and has a continuous derivative $p'_X(x)$ in some neighborhood of $x=0$, then 
\begin{equation}\label{eqa1}
\sup\nolimits_{x \in \mathbb{R}}\left| \mathbb{P}_{\theta}\left(2\,p_X(0) \sqrt{m}(M_m\,-\,\theta) \leq x\right) \,-\, 
\Phi(x)\right| \to 0 \quad \mbox{as} \quad m \to \infty,
\end{equation}
where $\Phi(x)$ is the standard Gaussian distribution function having density  $\varphi(x)$: 
\begin{equation}\label{N}
\Phi(x) = \int_{- \infty}^x \varphi(y) dy \quad \mbox{with} \quad \varphi(y)= \frac{1}{\sqrt{2 \pi}} \,\,e^{- y^2/2}.
\end{equation}

Instead of moment conditions now regularity assumptions on the density $p_X(x)$ are required:\\[0.5ex] 
{\bf Assumption A:} {\it The density $p_X(x)$ is symmetric around zero,
i.e., $p_X(- x) = p_X(x), x\in\mathbb{R}$ and $p_X(0) > 0$. Moreover, 
the density $p_X(x)$ has three continuous
bounded derivatives in some interval $(0, x_0), \,\, x_0 > 0$.}\\[0.5ex] 
Define \quad $p_0 = p_X(0) > 0,\quad p_1 = p'_X(0+) \quad \mbox{and} \quad  p_2 = p''_X(0+).$\\[0.5ex]
The regularity conditions in Assumption A are fulfilled, for example, for\\ 
$\bullet$ normal density~(\ref{N}),\\
$\bullet$ heavy tailed Student's $t$-distribution $S_{\nu}(x)$  with density function 
\begin{equation}\label{ST}
 s_{\nu}(x)= 
\frac{\Gamma((\nu+1)/2)}{\sqrt{\nu \pi} \, \Gamma(\nu/2)}\,\, 
\Big(1 + \frac{x^2}{\nu}\Big)^{- (\nu+1)/2},\quad \nu > 0, \quad x\in \mathbb{R},
\end{equation}
including Cauchy distribution in case $\nu=1$, where the degree of freedom parameter $\nu>0$ determines the heaviness  of the distribution tail,\\ 
$\bullet$ the triangular distribution with density
\begin{equation}\label{RD}
t_a(x) =  \frac{a - |x|}{a^2} {\rm \bf{1}}_{(-a\,,\,a)}(x),\,\,a>0
 \,\, \mbox{with} \,\, {\rm \bf{1}}_{A}(x):= 
\begin{cases} 1, & x\in A\\ 0, & x \notin A\end{cases},\,\,\,
A \subset \mathbb{R},
\end{equation}
$\bullet$ the  continuous uniform distribution or rectangular distribution with density
\begin{equation}\label{UD}
u_a(x) =  \frac{1}{2\,a} {\rm \bf{1}}_{(-a\,,\,a)}(x),\,\,a>0
\end{equation}
$\bullet$ and symmetric Laplace distribution  $L_{\mu}(x)$   having density 
\begin{equation}\label{eq6LA}
 l_{\mu}(x)=\frac{1}{\sqrt{2}\, \mu}\,e^{-\, \sqrt{2}\,|x|/ \mu}\,,\quad x\in \mathbb{R}, 
\quad \mu>0, \quad x \in \mathbb{R}.
\end{equation}\\[0.5ex]

The corresponding coefficients $p_0, p_1$, and $p_2$ in these examples are:
\begin{equation}\label{exam} \left.\hspace*{-0.2cm}\begin{array}{llll}
\bullet \,  \varphi(x): &\!\!\! p_0= 1/\sqrt{2\,\pi}, &\!\!\! p_1=0, &\!\!\! p_2= - 1/\sqrt{2\,\pi},\\
\bullet \,  s_{\nu}(x): &\!\!\! p_0= \frac{\Di \Gamma((\nu+1)/2)}{\Di \sqrt{v \pi}\,\Gamma(v/2)},\,\, &\!\!\! p_1=0,
&\!\!\! p_2= -\,\frac{\Di \Gamma((\nu+3)/2)}{\Di \sqrt{v \pi}\,\Gamma((v+2)/2)},\\
\bullet \,  t_a(x): &\!\!\! p_0 = a^{-1}, &\!\!\! p_1=- a^{-2},  &\!\!\! p_2= 0,\\
\bullet \,  u_a(x): &\!\!\! p_0 = (2\,a)^{-1}, &\!\!\! p_1= 0,  &\!\!\! p_2= 0,\\
\bullet \,  l_\mu(x): &\!\!\! p_0=1/(\sqrt{2} \, \mu), &\!\!\! p_1= - \mu^{- 2},\,\, &\!\!\! p_2=\sqrt{2} \mu^{- 3}.
\end{array}\right\} \end{equation}\\[0.5ex]

Under Assumption A \,\citet[Theorem~1]{Bu97} proved in relation (\ref{eqa1}) an asymptotic 
expansion in terms of orders $m^{-1/2}$ and $m^{-1}$ with  remainder $\mathcal{O}(m^{- 3/2})$ 
as $m \to \infty$. In the present paper we prove a similar second order expansion 
for the  sample median $M_{N_n}$ constructed from a sample with random sample size $N_n$. Therefore in Section 2 
we clarify the result of \citet{Bu97} in the sense that we get non-asymptotic relations for any integer $m \geq 1$ estimating
the closeness of the sample median $M_m$ and the corresponding second order expansion by inequalities. 
In Section 3 we give a transition proposition from non-random to random sample size and in 
Sections 4 and 5 the cases of Student $t$- and Laplace distributions 
as limit laws for the random median $M_{N_n}$ are considered. In Section~6 the Cornish-Fisher expansions for the quantiles of sample medians $M_{N_n}$ and $M_m$
are derived from the corresponding Edgeworth-type expansions.

\section{Non-Asymptotic Expansions for Sample Median }\label{s2} 
Let $[y]$ denote the integer part of value~$y$. Define 
\begin{equation}\label{eq20a}
m^* = 2\,[m/2] =  \left\{\begin{array}{ll}  m & \mbox{for even m,}\\
                                         m-1 \,\,& \mbox{for odd m.} \end{array}\right.
\end{equation}
\begin{proposition}\label{p1} Let Assumption A be satisfied, then for all $m \geq 2$: 
\begin{equation}\label{eq21} 
\sup\nolimits_{x \in \mathbb{R}}\left|\mathbb{P}_{\theta}(2 p_0 \sqrt{m^*} (M_m -\theta) \leq x ) - \Phi(x) - \frac{f_1(x)}{\sqrt{m^*}} 
- \frac{f_2(x)}{m^*} \right| \leq  \frac{C_1}{m^{3/2}}, 
\end{equation}
where $C_1$  does not depend on $m$,
\begin{equation}\label{eq33} 
f_1(x) = \frac{p_1 x |x| }{4 p_0^2}\varphi(x)\quad \mbox{and} \quad f_2(x) = \frac{x}{4}\Big(3 + 
x^2 +\frac{p_2 x^2}{6 p_0^3} - \frac{p_1^2 x^4}{8 p_0^4} \Big)\varphi(x).
\end{equation}
\end{proposition}

Since $0 < (m-1)^{- \alpha} - m^{- \alpha} \leq 2\,m^{- 3/2}$ for $m \geq 2$ and $\alpha= 1/2$ or $\alpha=1$ an immediate consequence of 
inequality (\ref{eq21}) is 
\begin{equation}\label{eq21c} 
\sup\nolimits_{x \in \mathbb{R}}\left|\mathbb{P}_{\theta}\Big(2 p_0 \sqrt{m^*} (M_m -\theta) \leq x \Big) - \Phi(x) - \frac{f_1(x)}{\sqrt{m}} 
- \frac{f_2(x)}{m} \right| \leq  \frac{C_2}{m^{3/2}},
\end{equation}
where  (\ref{eq21c})  for $m=1$ is trivial and $C_2$  does not dependent on $m$.\\[0.5ex]

\noindent
{\bf Remarks: 1.} If the parent distributions of the sample $\{X_1,...,X_m\}$ have the  normal density (\ref{N}),  
Student's $t$-density (\ref{ST}) or continuous uniform density (\ref{UD}), then with respect to 
(\ref{exam}) the first term $f_1(x)/\sqrt{m^*}$ vanishes since in these cases $p_1=0$. 
Therefore the convergence rate of the distribution of sample median $M_m$
to normality has order $m^{- 1}$. The triangular  density (\ref{RD}) and the Laplacian density (\ref{eq6LA}) have  discontinuous derivatives at $x=0$, nevertheless $p_1>0$ and the convergence rate to normality has the order $m^{- 1/2}$.

In \citet[Chapter 28.5]{Cr46} for asymptotic normality (\ref{eqa1}) it is required, that density $p_X(t)$ has a continuous derivative $p'_X(x)$ in some neighborhood of $x=0$.

{\bf 2.} As in \citet{Bu97} the natural normalizing factor in  (\ref{eq21}) is  $m^*$, i.e.,
$\sqrt{m-1}$ for odd $m\geq 3$  and  $\sqrt{m}$ for even $m$. He proved also for all $ m \geq 2 $  
$$ \left|\mathbb{P}_{\theta}(2 p_0 \sqrt{m^*} (M_{m^*} -\theta) \leq x)   - 
 \mathbb{P}_{\theta}(2 p_0 \sqrt{m^*} (M_{m^*+1} -\theta) \leq x )  \right| \leq  C m^{-3/2}.
$$
Hence, for the sample median $M_m$ each odd observation adds an
{\it{amount of information}} of order $m^{-3/2}$ and not $m^{-1}$ as usual if the normalizing factor by 
$(M_m - \theta)$ is $\sqrt{m}$. \\[0.5ex]

\noindent
{\it Proof of Proposition \ref{p1}:}
Following the detailed proof of \citet[Theorem~1]{Bu97} one has to change Stirling's formula 
of the Gamma functions $\Gamma(z)$ and $1/\Gamma(z)$ as $z \to \infty$  
 by inequalities, proved in \citet[Theorem~1.3]{Ne15}:
 \begin{equation}\label{Nem} \left. \begin{array}{lcl} \Gamma(z) &=& \sqrt{2\,\pi} \, z^{z - 1/2} \, e^{- z} \, 
(1 + \frac{\displaystyle 1}{\displaystyle 12 z} + \frac{\displaystyle 1}{\displaystyle 288 z^2} + R_{3}(z)),
 \\ [1ex]
\frac{\displaystyle 1}{\displaystyle \Gamma(z)} &=& \frac{\displaystyle 1}{\displaystyle \sqrt{2\,\pi}} \, 
z^{-z + 1/2} \, e^{z} \, (1 - \frac{\displaystyle 1}{\displaystyle 12 z} + 
\frac{\displaystyle 1}{\displaystyle 288 z^2} + \tilde{R}_{3}(z)),
\end{array}\right\} \quad z>0,
\end{equation}
with \,\,$\{|R_{3}(z)|,|\tilde{R}_{3}(z)| \leq c z^{- 3}$ \,\,and\,\, $ c= 
\frac{\displaystyle (1+ \zeta(3)) \Gamma(3)(2 \sqrt{3}+1)}{\displaystyle 2\, (2 \pi)^4} \leq 0.006$.\\
Here\,\, $\zeta(z)$ \,\,is the Riemann zeta function with\,\, $\zeta(3) \approx 1.202...$

Finally, when ever Taylor's formula 
is used with remainder in big $\mathcal{O}$ notation, then the remainder  has to be estimated in 
Lagrange form  by an inequality. The constants $C_1, C_2>0$ in (\ref{eq21}) and (\ref{eq21c}) depend only on $p_0, p_1, p_2$ 
and the upper bound of $p_X'''(x)$ in some interval $(0, x_0), \,\, x_0 > 0$.\hfill $\Box$

\section{Transfer Proposition from Non-Random to  Random Sample Sizes}
Suppose that  distribution functions of the random sample size $N_n$ satisfy the following  condition.\\[0.5ex]
{\bf Assumption B:} {\it There exist a distribution function $H(y)$ with $H(0+) = 0$,   
a function of bounded variation $h_2(y)$ with $h_2(0)=h_2(\infty)=0$, a sequence  $0<g_n \uparrow \infty$ and real 
numbers $b > 0$ and $C_3 > 0$   such~that for all $n \in \mathbb{N}$}
\begin{equation}\label{con2}\left. \begin{array}{ll}
\sup\nolimits_{y \geq 0} \left|\mathbb{P}\!\left(\!g_n^{-1} N_n \leq y\right) - H(y) \right| \leq C_3 n^{- b}, \quad&  0 < b \leq 1 \\[2ex]
\sup\nolimits_{y \geq 0} \left|\mathbb{P}\!\left(\!g_n^{-1} N_n \leq y\right) 
- H(y) - n^{-1} h_2(y) \right| \leq C_3 n^{- b}, \quad&  b > 1 
\end{array}\right\}
\end{equation}

\begin{theorem}\label{pro1} 
Let both Assumptions A and B be satisfied.  Then
the following inequality holds for all $n \in \mathbb{N}$ :
\begin{eqnarray}\label{p31}%
&&\sup\nolimits_{x \in \mathbb{R}}\Big| \mathbb{P}_{\theta}\Big(2 p_0 \sqrt{g_n \,N_n^*/N_n\,}\,\, 
(M_{N_n} - \theta) \leq  x\Big) - G_n(x,1/g_n)\Big| \nonumber\\
&&\hspace*{2cm} 
  \leq  C_2 \,\mathbb{E} \left(N_n^{- 3/2}\right)+  (C_3  D_n + C_4)\, n^{- b},
\end{eqnarray}
\begin{equation}\label{eqp1a}
G_n(x,1/g_n) =  \int\nolimits^{\infty}_{1/g_n}\Big(\Phi(x\sqrt{y})+\frac{f_1(x\sqrt{y})}{\sqrt{g_n y}} 
+\frac{f_2(x\sqrt{y})}{g_n y}\Big)d\Big(H(y) + \frac{h_2(y)}{n} \Big),
\end{equation}
\begin{equation}\label{p32} D_n = \sup\nolimits_x D_n(x) \leq D < \infty 
\end{equation}
$$
D_n(x)=  \int_{1/g_n}^\infty \left| \frac{\displaystyle \partial}{\displaystyle \partial y} \left( \Phi(x\sqrt{y}) 
+ \frac{f_1(x\sqrt{y})}{\sqrt{y g_n}} + \frac{f_2(x\sqrt{y})}{y g_n} \right) \right|dy ,
$$
where $f_1(z), f_2(z), h_{2}(y)$ are given in (\ref{eq33}) and  (\ref{con2}) and
\begin{equation}\label{Nn}
N_n^* = 2\,[N_n/2] = \left\{\begin{array}{ll}  N_n \quad & \mbox{for even realizations of $N_n$,}\\
                                         N_n - 1 \quad & \mbox{for odd realizations of $N_n$.} \end{array}\right. 
\end{equation}  
 The positive constants $C_2, C_3, C_4, D$ do not depend on $n$. 
\end{theorem}

\noindent    
{\bf Remarks: 1.} The  scaling factor $\sqrt{g_n\,N_n^*/ N_n}$ 
seems to be the  natural one in case of the median of a sample with a random sample size 
$N_n$ since the distribution of $N_n/g_n$ has a known limit distribution and $N_n^*$ the same structure as $m^*$ 
in \cite{Bu97}.

{\bf 2.} Without the quotient $N_n^*/N_n$ in the scaling factor $\sqrt{g_n\,N_n^*/ N_n}$ an additional term in the expansion occurs:
$$\frac{\Di 1}{\Di 2 n} \sum\nolimits_{u=1}^\infty \int\nolimits_{(2u-1)/g_n}^{2u/g_n} x \sqrt{y} \varphi(x \sqrt{y}) dH(y) 
\leq \frac{\Di 1}{\Di 2\sqrt{2 \pi e}\,n}.$$

{\bf 3.} The lower bound of the integral in (\ref{eqp1a}) depends on $g_n$ which can affect the coefficients at $1/\sqrt{g_n}$ and $1/g_n$ in the approximation. For example the proof of Theorem 4.2 in Section 4 shows that among other integrals
$$ g_n^{- 1} \int\nolimits^{\infty}_{1/g_n}\left|\frac{f_2(x\sqrt{y})}{y}\right| dH(y) \leq c g_n^{- b} \quad \mbox{if} \quad b< 1.$$
\\[1ex]

\noindent
{\it{Proof of Theorem \ref{pro1}:}}  The proof follows along the similar arguments 
of the more general transfer theorem in \citet[Theorem 3.1]{BGK13} under conditions of our
Theorem~\ref{pro1}. Then conditioning on $N_n$, we have 
\begin{eqnarray}\label{eqp1}
&&\hspace*{-0.2cm}\mathbb{P}_{\theta}\Big(2 p_0 \sqrt{g_n N_n^*/N_n}\big(M_{N_n} \!- \theta) \leq  x\Big)= 
\mathbb{P}_{\theta}\Big(2 p_0 \sqrt{N_n^*}(M_{N_n} \!- \theta) \leq  x \sqrt{N_n/g_n}\Big) \nonumber\\
&&\qquad \qquad  =\,\sum\nolimits_{m=1}^\infty \mathbb{P}_{\theta}\Big(2 p_0 \sqrt{m^*}(M_{m} - \theta) \leq  x \sqrt{m/g_n}\Big)\, \mathbb{P}(N_n = m).
\end{eqnarray}
Using now  (\ref{eq21c}) with $\Phi_m(z):= \Phi(z) + m^{- 1/2} f_1(z) + m^{-1} f_2(z)$:\\
$$\begin{array}{l}\label{eqp3a}
{\Di \sup\nolimits_x\sum\nolimits_{m=1}^\infty \left|\mathbb{P}_{\theta}\Big(2 p_0 \sqrt{m^*}(M_{m} - \theta) \leq x\sqrt{m/g_n} \Big)\,- \Phi_m(x\sqrt{m/g_n})\right| \,\mathbb{P}(N_n = m)}\nonumber\\
 {\Di \qquad \qquad \qquad \stackrel{(\ref{eq21c})}{\leq} C_2\,\sum\nolimits_{m=1}^\infty  m^{- 3/2} \,\mathbb{P}(N_n = m)
=   C_2  \, \mathbb{E}(N_n^{- 3/2}). }
\end{array}$$
Taking in account  $\mathbb{P}\Big(N_n/g_n < 1/g_n\Big) = \mathbb{P}\Big(N_n < 1\Big) = 0$ we obtain
\begin{eqnarray}\label{eqp3}
&&\sum\nolimits_{m=1}^\infty \Phi_m(x\sqrt{m/g_n})\mathbb{P}(N_n = m)
= \mathbb{E}_{\theta} \left(\Phi_{N_n}(x\sqrt{N_n/g_n})\right) \nonumber\\
&&\hspace*{1.5cm}=\int\nolimits_{1/g_n}^\infty \Delta_n(x,y) d\mathbb{P}(N_n/g_n \leq y) = G_n(x,1/g_n) + I,\nonumber
\end{eqnarray}
where \quad $\Delta_n(x,y) :=  \Phi(x \sqrt{y}) + f_1(x\sqrt{y})/\sqrt{g_n y} +f_2(x\sqrt{y})/(g_n y)$,  
$G_n(x,1/g_n)$ is defined in (\ref{eqp1a}) and
$$
I= \int_{1/g_n}^\infty  \Delta_n(x,y)d\Big(\mathbb{P}(N_n/g_n \leq y) - H(y)-\frac{h_2(y)}{n}\Big).
$$
Estimating integral $I$ we use integration by parts for Lebesgue-Stieltjes integrals.
\begin{eqnarray*}
|I| &\leq & \sup\nolimits_x \left. \lim_{L\to\infty}|\Delta_n(x,y)|\, 
\big|\mathbb{P}\big(N_n/g_n \leq y\big) - H(y) - n^{-1} h_2(y)\big|\right|_{y=1/g_n}^{y=L}\nonumber\\ 
&&+ \, \sup\nolimits_x \int_{1/g_n}^\infty \Big|\frac{\Di \partial}{\Di \partial\,y}\Delta_n(x,y)\Big|\left|\mathbb{P}\big(N_n/g_n \leq y\big)
- H(y) -n^{-1}h_2(y)\right|dy. 
\end{eqnarray*}
First we calculate  $(\partial/\partial\,y)\Delta_n(x,y)$. Obviously $\frac{\Di \partial}{\Di \partial y} \Phi(x\,\sqrt{y})=  \, \frac{\Di x}{\Di 2}\, y^{- 1/2} \varphi(x\,\sqrt{y})$ and
\begin{equation}\label{del1}
\frac{\partial}{\partial\,y}\Big(\frac{f_1(x\sqrt{y})}{\sqrt{y} }\Big) = \frac{{\rm{sgn}}(x) \, q_1(x\sqrt{y})}{2\,y^{3/2}}\,, \qquad \frac{\partial}{\partial\,y}\Big(\frac{f_2(x\sqrt{y})}{y}\Big) = \frac{q_2(x\sqrt{y})}{8 \, y^{2}},
\end{equation}
\begin{equation*}
q_1(z)=a_0 (1 - z^2 )z^2\,\varphi(z), \quad q_2(z) = \Big(a_2 z^6 - (a_1 + 3 a_2) z^4 + (a_1 - 3)z^2 - 3\Big) z \varphi(z),
\end{equation*}
where $a_0=p_1/(4 p_0)$, \,$a_1=1+p_2/(6 p_0^3)$ \,\, and \,\, $a_2=p_1^2/(8 p_0^4)$,\, see (\ref{eq33}).

The functions $f_k(z)$ and $q_k(z)$, $k=1, 2$, are bounded, we suppose
\begin{equation}\label{fq}
\sup\nolimits_z |f_k(z)| \leq c_k^* < \infty \quad \mbox{and} \quad \sup\nolimits_z |q_k(z)| \leq c_k^{**} < \infty,\quad k=1,2. 
\end{equation}

To estimate $D_n$ defined in \eqref{p32} we consider $D_n(x)$ for $x \not= 0$ since $D_n(0)=~0$. Because $0 \leq \int_{1/g_n}^\infty (\partial/\partial\,y)\Phi(x \sqrt{y})dy = 1 - \Phi(x/\sqrt{g_n}) \leq 1/2$ \,\,for\,\, $x > 0$ and 
$\int_{1/g_n}^\infty |(\partial/\partial\,y)\Phi(x \sqrt{y})|dy = \Phi(x/\sqrt{g_n}) \leq 1/2$ \,\,for \,\, $x < 0$  we find
with~\eqref{fq} \,\,   
$D_n(x) \leq 1/2 +c_1^{**} + c_2^{**}/8 $ for $x \not= 0$. Therefore  inequality (\ref{p32}) holds with $D = 1/2+ c_1^{**} + c_2^{**}/8$. It follows now from \eqref{con2} and \eqref{fq} that
\begin{equation*}
|I| \leq    (1+ c_1^* + c_2^*)\, C_3 \, n^{- b}   + (1/2+ c_1^{**} + c_2^{**}/8)\Big) \, C_3 \, n^{- b} 
\end{equation*}
and $C_4 =(1+ c_1^* + c_2^*)\, C_3$. Theorem \ref{pro1} is proved.
\hfill $\Box$\\[1ex]

\begin{theorem}\label{pro2} 
Under the conditions of Theorem \ref{pro1}  
and the additional conditions  to functions $H(.)$ and $h_2(.)$, depending on the convergence rate $b>0$ in (\ref{con2}):
\begin{equation}\label{add2a}\left. \begin{array}{llll}
i: \quad & H(1/g_n) \leq c_1\,g_n^{- b} \qquad &  \mbox{for}  & b > 0,\\[0.5ex]
ii: \quad &  \int\nolimits_0^{1/g_n} y^{-\,1/2}  dH(y) \leq c_2\,g_n^{- b + 1/2} \qquad &  \mbox{for}  & b > 1/2,\\[0.5ex]
iii: \quad & \int\nolimits_0^{1/g_n} y^{-\,1}  dH(y)\leq c_3\,g_n^{- b + 1} \qquad & \mbox{for} & b>1, 
\end{array}\,\right\}\end{equation}
\begin{equation}\label{add2b}\left. \begin{array}{llll}
i: & h_2(0) = 0, \quad \mbox{and} \quad |h_2(1/g_n)| \leq c_4\,n\,g_n^{- b}& \mbox{for} & b >1,\\[0.5ex]
ii: & \int\nolimits_0^{1/g_n} y^{-\,1} |h_2(y)| dy \leq c_5\,n \,g_n^{- b}& \mbox{for} & b>1,
\end{array}\,\right\}\end{equation}
we obtain for the function $G_n(x,1/g_n)$ defined in (\ref{eqp1a}): 
\begin{equation}\label{eqpb}
\sup\nolimits_x \big|G_n(x,1/g_n) - G_{n,2}(x)\big| \leq C \, g_n^{- b} + \sup\nolimits_x \left(|I_1(x,n)| + |I_2(x,n)|\right)
\end{equation}
with
\begin{equation}\label{gn}
G_{n,2}(x) = \left\{ \begin{array}{lr} \int\limits^{\infty}_{0} \Phi(x\sqrt{y}) d H(y), & \hspace*{-1cm} 0<b \leq 1/2, \\
 \int\limits_{0}^{\infty} \Big(\!\Phi(x\sqrt{y}) + \frac{\Di f_1(x\sqrt{y})}{\Di \sqrt{g_n y}} \Big)
 dH(y) =: G_{n,1}(x), & \hspace*{-1cm} 1/2 < b \leq 1,\\
G_{n,1}(x) + \int\limits^{\infty}_{0} \frac{\Di f_2(x\sqrt{y})}{\Di g_n y}dH(y) 
 + \int\limits_{0}^{\infty} \frac{\Di \Phi(x\sqrt{y})}{\Di n}  dh_2(y), &\hspace*{-0.9cm} b > 1,\end{array} \right.
\end{equation}
\begin{equation}\label{gnI1}
I_1(x,n) = \int\nolimits_{1/g_n}^\infty  \Big(\frac{f_1(x\sqrt{y})\,\,{\rm \bf{1}}_{(0 , 1/2]}(b)}{\sqrt{g_n y\,}} 
+ \frac{f_2(x \sqrt{y\,})}{g_n\,y}\Big)  dH(y) \quad \mbox{for} \quad b\leq 1
\end{equation}
and 
\begin{equation}\label{gnI1b} 
I_2(x,n) = \int\nolimits_{1/g_n}^\infty \Big(\frac{f_1(x\sqrt{y})}{n\,\sqrt{g_n y\,}} 
+ \frac{f_2(x\sqrt{y})}{n\,g_n y}\Big)  dh_2(y) \quad \mbox{for} \quad b > 1.
\end{equation}
\end{theorem}

\noindent 
{\bf Remarks:} If $b >1/2$ then  (\ref{add2a}ii) implies (\ref{add2a}i). If $b >1$ then  (\ref{add2a}iii) implies (\ref{add2a}ii) and (\ref{add2a}i). 
 Conditions (\ref{add2a}) and (\ref{add2b}) lead to the range of the integrals in (\ref{gn}) which  ensures~(\ref{eqpb}).
The length of the asymptotic expansion is defined by (\ref{gn}).\\[2ex]

\noindent 
{\it Proof of Theorem \ref{pro2}:} Using  condition (\ref{add2a}i) we find
$$\int\nolimits^{1/g_n}_{0} \Phi(x\sqrt{y}) d H(y) \leq \int\nolimits^{1/g_n}_{0}  d H(y)
 = H(1/\sqrt{g_n}) - H(0) \leq c_1g_n^{-b}.$$
It follows from \eqref{fq},  (\ref{add2a}ii) and (\ref{add2a}iii) that for $k=1, 2$
$$\int\nolimits^{1/g_n}_{0} |f_k(x\sqrt{y})| y^{-k/2}d H(y) \leq c_k^*
\int\nolimits^{1/g_n}_{0}  y^{-k/2}d H(y) \leq c_k^*\, c_{k+1}(r) g_n^{- b + k/2}.$$
Integration by parts, $|z|\varphi(z)/2 \leq c^*=(8\,\pi \, e)^{-1/2}$, (\ref{add2b}i) and (\ref{add2b}ii) lead to
$$\left|\int\nolimits^{1/g_n}_{0} \!\Phi(x\sqrt{y})  dh_2(y)\right| \leq  |h_2(1/g_n)| + c^*
\int\nolimits^{1/g_n}_{0} \! y^{-1}  |h_2(y)| dy \leq (c_4 +c^*c_5)n\,g_n^{- b}. $$
Taking into account  (\ref{eqp1a}), (\ref{gn}), (\ref{gnI1}) and (\ref{gnI1b}) we obtain (\ref{eqpb}).   \hfill $\Box$\\[1ex]

In the next two sections we use both Theorems 3.1 and 3.2 when the scale mixture
 $G(x) = \int_0^\infty \Phi(x\,\sqrt{y}) d H(y)$ as limiting distribution of $M_{N_n}$ can be expressed 
in terms of the well-known distributions. We obtain non-asymptotic results like in Proposition \ref{p1} 
for the sample median $M_{N_n}$, using second order approximations for both
the statistic $M_m$ and for the random sample size $N_n$. In both cases the jumps of the distribution function of the random sample size $N_n$ only affect the function $h_2(y)$ in formula (\ref{con2}).

\section{Student's~Distribution as Limit for Random Sample Median $M_{N_n}$}
Let the sample size $N_n(r)$ be the negative binomial distributed (shifted by~1)   with parameters $1/n$ and $r > 0$,  having probability mass function
\begin{equation}\label{eq5} 
\mathbb{P}(N_n(r)=j)= \frac{\Gamma(j+r-1)}{\Gamma(j) \, \Gamma(r)}\left(\frac{1}{n} \right)^{r}\left( 1 - \frac{1}{n} \right)^{j-1}, \,\,\,  j =1, 2, ... 
\end{equation}
with $g_n=\mathbb{E}(N_n(r))= r\,(n-1) +1$. \citet[Section 2.1]{ST16} pointed out that the negative binomial distribution is one of the two leading cases for count models, it accommodates the over-dispersion typically observed in count data (which the Poisson model cannot) and they showed in a general unifying framework
\begin{equation}\label{Nr}                              
\lim\nolimits_{n\to\infty} \sup\nolimits_x\left|\mathbb{P}(N_n(r)/g_n \leq x) - G_{r,r}(x)\right| = 0,
\end{equation} 
where $G_{r,r}(x)$ is the Gamma distribution function   with the shape parameter which coincides with the scale parameter and equals $r >0$, having density
\begin{equation}\label{Gam}
g_{r,r}(x)=\frac{r^r}{\Gamma(r)} x^{r - 1} e^{- r x}\,{\rm \bf{1}}_{(0\,,\,\infty)}(x),\quad x \in \mathbb{R}.
\end{equation}
The statement (\ref{Nr}) was proved earlier in \citet[Lemma 2.2]{BK05}.

The convergence rate in (\ref{Nr}) for $r > 0$ is given in \citet[Formula (21)]{BGK13} or \citet[Formula (17)]{GZK17}:
\begin{equation}\label{Nr1} 
 \sup\nolimits_x\left|\mathbb{P}(N_n(r)/g_n \leq x) - G_{r,r}(x)\right| \leq C_r n^{- \min\{r, 1\}}.
\end{equation}
In \citet{ST16} and \citet{GZK17} the negative binomial random variable $\tilde{N}_n(r)$ is not shifted: $\tilde{N}_n(r) = N_n(r) -1 \in \{0, 1, 2, ...\}$ with $\mathbb{E}\tilde{N}_n(r) = r(n-1)$. Then we have  $\mathbb{P}(\tilde{N}_n(r) \leq 0) - G_{r,r}(0)= n^{- r} \to 0$ as $n\to \infty$ instead of $\mathbb{P}(N_n(r) \leq 0) - G_{r,r}(0)=0$. Moreover
$$ \mathbb{P}\left(\frac{\tilde{N}_n(r))}{r(n-1)} \leq x\right) = \mathbb{P}\left(\frac{N_n(r)}{g_n} \leq x +\frac{1-x}{g_n}\right).$$ The statements (\ref{Nr}) and (\ref{Nr1}) still hold when $\tilde{N}_n(r)$ is shifted by a fixed integer. 
From Taylor expansion with Lagrange remainder term it follows that for $r>1$ 
$$ \left|G_{r,r}\left(x +\frac{1-x}{g_n}\right) - G_{r,r}(x) - g_{r,r}(x)\frac{1-x}{g_n}\right| \leq C \,g_n^{- \min\{r , 2\}}.
$$
Hence, for $r>1$ shifting  $\tilde{N}_n(r)$ has influence of a term by $g_n^{-1}$.  Second order asymptotic expansions for $N_n(r)$ where proved  in \citet[Theorem~1]{CMU20}:
\begin{proposition}\label{p2} Let $r > 0$,  discrete random variable $N_n(r)$ have probability mass function 
(\ref{eq5}) and $g_n:= \mathbb{E}N_n(r) = r (n-1) +1$. For $x>0$ and all $n \in \mathbb{N}$  there exists a real number $C_3(r) > 0$ such that
\begin{equation}\label{eq9a} 
\sup\nolimits_{x \geq 0} \left|\mathbb{P} \left(\frac{N_n(r)}{g_n} \leq x\right) - G_{r,r}(x) - \frac{h_{2;r}(x)}{n} \right| \leq C_3(r) \,n^{- \min\{r , 2\}},
\end{equation} 
where
\begin{equation}\label{eq9c}
h_{2;r}(x)  = \begin{cases}  0 & for \,\, r \leq 1,\\
\frac{\Di g_{r,r}(x) \left( (x-1)(2 - r) + 2 Q_1\big(g_n \,x\big)\right)}{\Di 2\,r}& for \,\, r>1,
\end{cases}
\end{equation}
\begin{equation}\label{eq9aa} 
Q_1(y) = 1/2 - (y - [y]) \,\, \mbox{and}\,  [.] \, \mbox{denotes the integer part of a number}. 
\end{equation}
\end{proposition}

\noindent
{\bf Remark:}  The jumps of the sample size $N_n(r)$ have an effect only on the function $Q_1(.)$ in the term $h_{2;r}(x)$. The function $Q_1(y)$ is periodic with period 1,
it is right-continuous with jump height 1 at each integer point $y$. The Fourier series expansion of $Q_1(y)$ at all non-integer points $y$ is
\begin{equation}\label{eqA6}
 Q_1(y) = 1/2 -(y - [y]) = \sum\nolimits_{k=1}^{\infty} \frac{\sin(2\,\pi\, k \, y)}{k\,\pi}\,,\quad y\not=[y],
\end{equation}
see formula 5.4.2.9 in \citet[ p.~726]{PBM92} with $a=0$.\\[1ex]

In Theorem \ref{pro1} an estimate for the negative moment $\mathbb{E}(N_n)^{-3/2}$ of the random sample size $N_n$ is required.
Proposition \ref{p2} is used in \citet[Corollary~2]{Ben20} to obtain an asymptotic expansion of negative moments $\mathbb{E}(N_n(r))^{-p}$ for $r >1$ and $0<p \leq r-1$. Such expansions are applied in the mentioned paper to to analyze asymptotic deficiencies 
 and risk functions of estimates based on random-size samples. An improved result is given here, i ncluding the correct bounds for $0<p \leq r-1$:
\begin{corollary}\label{cor42}
Let  $r >0$ and $p>0$. Then for all $n\geq 2$ the following  expansions hold for  negative moments:
\begin{equation}\label{c4}
\mathbb{E}(N_n(r))^{-p} = {\Di \begin{cases} {\Di \frac{\Di r^p\,\Gamma(r-p) }{\Gamma(r)\,g_n^p} - \frac{(2-r)\,p\,r^p (p+1)\, \Gamma(r-p-1)}{2\,r\,\Gamma(r)\,n\,g_n^p }} + R_{1;n}^*, & p < r -1,\\[1.5ex]
 {\Di \frac{\Di r^p\,\Gamma(r-p) }{\Gamma(r)\,g_n^p} - \frac{(2-r)\,p\,r^p \,\ln (g_n/r)}{2\,r\,\Gamma(r)\,n\,g_n^p }} \,
 + \, R_{2;n}^*, & p = r-1,  \\[1.5ex]
 {\Di \frac{\Di r^p\,\Gamma(r-p) }{\Gamma(r)\,g_n^p} \,  + \, R_{3;n}^*}, &\hspace*{-0.6cm} r-1 < p < r, \\[1.5ex]
{\Di \frac{\Di r^r\,\ln (g_n/r)}{ \Gamma(r)\,g_n^r}\,+ \, R_{4;n}^*,} & p=r, \\[1.5ex]
{\Di \frac{r^r}{\Gamma(r)\,(p-r)\,g_n^r} \,+\, R_{5;n}^*,}  & p > r,
\end{cases}} \end{equation}  
where $|R_{k;n}^*| \leq c_k^*(p,r) \, g_n^{- \min\{r,2\}}$ for some constants $c_k^*(p,r) < \infty$, \,\,$k=1, 2, ...,5$.\\[1ex]
\end{corollary} 

\noindent
{\bf Remark:} The leading terms in \eqref{c4} and the bound \eqref{eq9a} lead to the estimate
\begin{equation}\label{NM2a} \mathbb{E}\big(N_n(r)\big)^{- p} \leq C(r,p) 
\begin{cases} n^{- \min\{r, p, 2\}}, &  p \not= \min\{r , 2\}\\ \ln(n) \, n^{- \min\{r, p, 2\}}, &   p  = \min\{r , 2\} 
 \end{cases} 
\end{equation}

\noindent
{\it Proof of Corollary~\ref{cor42}:} Integrating by parts and substituting $y/g_n = x$ , we obtain 
\begin{eqnarray}\label{cor1}
\mathbb{E}(N_n(r))^{-p} = \int\nolimits_1^\infty \frac{1}{y^p} d\mathbb{P}\big(N_n(r) < y\big) = \frac{p}{g_n^p} 
\int\nolimits_{1/g_n}^\infty \frac{1}{x^{p+1}} \mathbb{P}\left(\frac{N_n(r)}{g_n} < x\right)dx&& \nonumber\\
= \frac{p}{g_n^p}  \int\nolimits_{1/g_n}^\infty \frac{1}{x^{p+1}} \left(G_{r,r}(x) + \frac{h_{2;r}(x)}{n}\right)dx  +R_1(n) =I_1 + I_2 + R_1(n), \qquad &&
\end{eqnarray}
where with \eqref{eq9a}  of Proposition \ref{p2} 
\begin{equation*}
|R_1(n)| \leq \frac{p}{g_n^p}  \int\nolimits_{1/g_n}^\infty \frac{1}{x^{p+1}} \left|\mathbb{P}\left(\frac{N_n(r)}{g_n} < x\right) - G_{r,r}(x) - \frac{h_{2;r}(x)}{n}\right|dx \leq \frac{C_3(r)}{n^{\min\{r , 2\}}}.
\end{equation*}
Next we calculate the first part $I_1$ of the integral in \eqref{cor1}:
\begin{equation*}
 I_1 =\frac{p}{g_n^p} \,\int\nolimits_{1/g_n}^\infty \frac{G_{r,r}(x)}{x^{p+1}}dx = I_{1;p}(n) + R_2(n)
\end{equation*}
with $R_2(n) = G_{r,r}(1/g_n) \leq r^r g_n^{-r}/\Gamma(r+1)$
\begin{equation}\label{cor2}
I_{1;p}(n) =  \frac{r^r}{\Gamma(r)\,g_n^p}  \int\nolimits_{1/g_n}^\infty  \frac{e^{- r\,x}}{x^{p+1-r}}  dx = \begin{cases}  \frac{\Di r^p\,\Gamma(r-p) }{\Di \Gamma(r)\,g_n^p}+ R_{3;p}(n), & p < r,\\[1.5ex]
\frac{\Di r^r\,\ln (g_n/r)}{\Di \Gamma(r)\,g_n^r} + R_{4;p}(n), & p=r, \\[1.5ex]
\frac{\Di r^r}{\Di \Gamma(r)\,(p-r)\,g_n^r} + R_{5;p}(n),  & p > r,
\end{cases}
\end{equation}
where for $p < r$ 
\begin{equation}\label{corNB}
R_{3;p}(n) =  \frac{\Di r^r}{\Di \Gamma(r)\,g_n^{p}} {\Di \int\nolimits_0^{1/g_n}} x^{r-p-1} e^{- r\,x}  dx  \leq \frac{\Di r^r}{\Di \Gamma(r)\,(r - p)\,g_n^{r}}, \qquad  p < r.\\[1.5ex]
\end{equation}
In case $p \geq r $ we split the integral in $I_1(n,p)$ into three parts, the first one leads to the leading term in \eqref{cor2}   for $p=r$ and $p>r$, respectively: 
\begin{equation*} 
\int\nolimits_{1/g_n}^\infty \frac{e^{- r\,x}}{x^{p+1-r}} dx=
{\Di \int\nolimits_{1/g_n}^{1/r}} x^{r-p-1}   dx -  {\Di \int\nolimits_{1/g_n}^{1/r}}  \frac{1 - e^{- r\,x} }{x^{p+1-r}}  dx + {\Di \int\nolimits_{1/r}^\infty} \frac{e^{- r\,x}}{x^{p+1-r}}  dx.
\end{equation*}
Then we obtain
\begin{equation*} \hspace*{-0.2cm}\begin{array}{rr} 
R_4(n,p) =  {\Di \frac{r^r}{\Gamma(r)\,g_n^{p}} \left( \int\nolimits_{1/g_n}^{1/r} \frac{1 - e^{- r\,x}}{x}  dx +  \int\nolimits_{1/r}^\infty  \frac{ e^{- r\,x}}{x}  dx \right)
	\leq \frac{r^r\,(1 +  e^{-1})}{\Gamma(r)\,g_n^{r}},} & p=r, \\[1.5ex]
|R_5(n,p)|  = {\Di \frac{r^r}{\Gamma(r)\,g_n^{p}} \left|-\, \frac{r^{p-r}}{p-r}\,-\, \int\nolimits_{1/g_n}^{1/r} \frac{1 - e^{- r\,x}}{x^{p-r+1}}dx + \int\nolimits_{1/r}^\infty \frac{e^{- r\,x}}{x^{p-r+1}} dx\right|}\quad& \\[1.5ex]
 \leq  {\Di \frac{r^{p}}{\Gamma(r)\,g_n^{p}}  \left(\frac{1}{p-r} + \frac{{\rm \bf{1}}_{(r , r+1)}(p)}{r + 1 - p}  + \ln(g_n/r)\,{\rm \bf{1}}_{\{p=r+1\}}(p) \,+\, e^{-1} \right)}\quad& \\
\qquad \quad + \, \frac{\Di r^{r+1}}{\Di \Gamma(r)\,(p - r - 1)\,g_n^{r+1}}  \,{\rm \bf{1}}_{\{p>r+1\}}(p), \qquad\quad&  p> r.\end{array}
 \end{equation*}
Now we calculate the second part $I_2$ of the integral in \eqref{cor1} in case of $r>1$:
\begin{equation}\label{cor2a}
I_2 =\frac{p}{g_n^p\,n}\int\nolimits_{1/g_n}^\infty  \frac{h_{2;r}(x)}{x^{p+1}}dx  = \frac{p\,r^r\,(2 - r)}{2\,r\,\Gamma(r)\,g_n^p\,n}  \int\nolimits_{1/g_n}^\infty \frac{e^{-r\,x}}{x^{p + 2 - r}} (x-1)dx + I_{2,p}(n).
\end{equation}
First we show that the integral $I_{2,p}(n)$ has the order of the remainder:
\begin{equation}\label{cor3}
|I_{2,p}(n)| = \frac{p}{r\,g_n^p\,n}\left|\int\nolimits_{1/g_n}^\infty \frac{e^{-r\,x}}{x^{p + 2 - r}} Q_1\big(g_n \,x\big)dx\right| \leq c(p,r) g_n^{- r}. 
\end{equation}
Let  $0<p < r-1$ where $r >1$. Then
$$ I_{2,p}(n) = \frac{p}{r\,g_n^p\,n}\int\nolimits_{0}^\infty \frac{e^{-r\,x}}{x^{p + 2 - r}} Q_1\big(g_n \,x\big)dx + R_6(n)$$
where since $ g_n \leq r\, n$ for $r>1$
$$|R_6(n)|= \frac{p}{r\,g_n^p\,n}\left|\int\nolimits_{0}^{1/g_n} \!\! \frac{e^{-r\,x}}{x^{p + 2 - r}} Q_1(g_n \,x)dx \right|\leq  \frac{p}{2 r\,g_n^p\,n} \int\nolimits_{0}^{1/g_n} \!\!\!\!\frac{dx}{x^{p + 2 - r}}  = \frac{g_n^{- r}}{2(r-p-1)}. $$
and considering \eqref{eqA6} and interchange integral and sum 
$$J_{2;p}^*(n) = \int\nolimits_{0}^\infty \frac{e^{-r\,x}}{x^{p + 2 - r}} Q_1\big(g_n \,x\big)dx = \sum\nolimits_{k=1}^{\infty} \frac{1}{k\,\pi}\int\nolimits_{0}^\infty \frac{e^{-r\,x}}{x^{p + 2 - r}}\,\sin(2\,\pi\, k \, g_n\, x) dx.$$
Applying formula 2.5.31.4 in \citet[p.~446]{PBM92} with $\alpha = r - p-1,\, p=r$ and $b=2\pi k g_n$ then 
\begin{eqnarray*}
J_{2;p}^*(n) &=&\sum\nolimits^{\infty}_{k=1}\frac{1}{k\,\pi}\int\nolimits_{0}^\infty x^{ r -p -2}e^{-r\,x}\,\sin(2\,\pi\, k \, g_n\, x) dx \\
&=&\frac{\Gamma(r - p -1)}{\pi}\,\sum\nolimits^{\infty}_{k=1}\frac{\sin\big((r - p - 1)\arctan\left(2 \pi k g_n/r\right)\big)}{k\,\big(\big(2\pi k g_n\big)^2+ r^2 \big)^{(r - p -1)/2}}.
\end{eqnarray*}
Hence 
$$|J_{2;p}^*(n)| \leq \frac{ \Gamma(r - p - 1)}{\pi} \sum\nolimits^{\infty}_{k=1} \frac{1}{k\,(2\pi\,k \,g_n)^{r-p-1}} \leq 
\frac{ \Gamma(r - p - 1)}{\pi (2\pi)^{r-p-1}}\, \zeta(r-p)\,g_n^{- r+p+1}$$
with Riemann zeta function $\zeta(.)$ and $\frac{\Di 1}{\Di r-p-1}\leq\zeta(r-p)\leq \frac{\Di r-p}{\Di r-p-1}< \infty$ and 
$$ |I_{2,p}(n)| \leq \left(\frac{p\, \Gamma(r - p - 1)\,\zeta(r-p)}{\pi (2\pi)^{r-p-1}} + \frac{r}{r-p-1} \right)\,g_n^{- r}, \qquad 0<p<r-1. $$

In case $0<p=r-1$ the Fourier series expansion \eqref{eqA6} of $Q_1(y )$ and integration by parts lead to
\begin{equation*}
I_{2,p}(n) = \frac{p}{r\,g_n^p\,n}\,\sum\limits^{\infty}_{k=1}\frac{1}{k^2\,2\,\pi^2\,g_n}\left(\frac{g_n}{e^{r/g_n}} -  \int\nolimits_{1/g_n}^\infty \left(\frac{1}{x^2} + \frac{1}{x}\right)e^{-r\,x}\cos(2\,\pi\, k \, g_n\, x) dx\right) 
\end{equation*}
and
\begin{equation*}
|I_{2,p}(n) |\leq
\frac{p}{r\,g_n^p\,n}\,\sum\nolimits^{\infty}_{k=1}\frac{2\, g_n + e^{-1} g_n}{k^2\,2\,\pi^2\,g_n} \leq
\frac{(2+e^{-1}) p}{2 \pi^2} \, \zeta(2)\, g_n^{- r}, \quad p=r-1.
\end{equation*}

If $p>r-1$ using $|Q_1(y) \leq 1/2$ we find
\begin{equation*}
|I_{2,p}(n) |\leq \frac{p}{2\,r\,g_n^p\,n} \int\nolimits_{1/g_n}^\infty x^{-r + p+1}  dx \leq  \frac{p}{2\,(p+1-r)} g_n^{-r}, \quad p > r-1,
\end{equation*}
and \eqref{cor3} is proved.

It remains to calculate the first term on the right-hand side of \eqref{cor2a}, say  $I_{3,p}(n) =I_2 - I_{2,p}(n)$. Since the integrals in $I_1$ and $I_{3,p}(n)$ have the same structure, one get with the above method
\begin{equation}\label{cor4}
I_{3,p}(n) = \begin{cases} {\Di \frac{p\,r^p\,(2-r)}{2\,r\,\Gamma(r) g_n^p\,n}} (\Gamma(r-p) - r\,\Gamma(r-p-1)) + R_7(n), & p < r-1, \\[1.5ex]
{\Di \frac{p\,r^{r}\,(2-r)}{2\,r\,\Gamma(r) g_n^{r-1}\,n}} (- r\,\ln(g_n/r))+R_8(n), & p= r - 1,\\[1.5ex]
R_9(n), & p > r - 1,
\end{cases}
\end{equation}
where $|R_k(n)| \leq c_k(r,p) g_n^{- r}$, with some constants $c_k(r,p)$,  $k=7, 8, 9$.

Estimates \eqref{cor1}, \eqref{cor2} and \eqref{cor4} lead to \eqref{c4} and Corollary \ref{cor42} is proved. \hfill $\Box$\\[1ex]

If the statistic $T_m$ is asymptotically normal the limit distribution of the standardized statistic $T_{N_n(r)}$ with random size $N_n(r)$ is Student's $t$-distri\~bution  $S_{2 r}(x)$  having density (\ref{ST}) with $\nu=2 r$,
see \citet{BK05} or \citet{ST16}. 

\noindent 
\begin{theorem}\label{th2}
Let $r>0$. Consider the sample median $M_{N_n}$ with random sample size $N_n = N_n(r)$ having probability 
mass function (\ref{eq5}) and $g_n=\mathbb{E}N_n(r) = r (n-1) +1$. 
If inequalities (\ref{eq21c}) and (\ref{eq9a}) hold for the mean $M_m(X_1,...,X_m)$ and the random sample size $N_n(r)$, respectively,  then there exists a constant $C_r$ such that 
\begin{equation}\label{teq2}
\Big|\mathbb{P}\Big(2 p_0 \sqrt{\frac{g_n N_n^*}{N_n}}(M_{N_n} - \theta)\leq  x\Big) - S_{2r;2}(x;n))\Big|\! \leq\! 
\begin{cases} C_r n^{- \min\{r , 3/2\}}\!, & \!\!\!r\not=3/2 \\
C_r \ln n \,\,n^{-3/2}\!, & \!\!\!r=3/2 
\end{cases}
\end{equation} 
for all $n \in \mathbb{N}$ uniformly in $x\in \mathbb{R}$, where $N_n^*$ is defined in (\ref{Nn}),
\begin{equation}\label{teq3a}
S_{2r;2}(x;n) = S_{2r}(x) + s_{2r}(x)\,\left(A_{1;r}(x) g_n^{- 1/2} \,+ A_{2;r}(x) g_n^{- 1} \right),
\end{equation}
$$ A_{1;r}(x) = \frac{p_1\,x\ |x|}{4\,p_0^2} \,{\rm \bf{1}}_{(1/2 , \infty)}(r)\qquad and \qquad$$
$$A_{2;r}(x) = \frac{x}{4}\left(\frac{(5 - r) x^2 +5 r + 2}{2 r - 1} + 
x^2 \Big(1 + \frac{p_2}{6 p_0^3}\Big) - \frac{x^4  p_1^2 (2 r + 1)}{8 p_0^4 (2 r + x^2)}\right){\rm \bf{1}}_{(1 , \infty)}(r).$$
\end{theorem}
{\bf Remark:} Under the condition (\ref{Nr1}) with $r>1/2$ a first order expansions of 
$\mathbb{P}_{\theta}\big(2 p_0 \sqrt{g_n}(M_{N_n} - \theta) \leq  x\big)$ was announced in the conference paper \citet{BKZ16}. 
Note that the convergence rate in Theorems 3.1 and 3.2 as well as in Corollaries 3.1 and 3.2  in case $1/2 < r < 1$ has to be $\mathcal{O}(n^{-r})$ instead of $\mathcal{O}(n^{- 1})$ as announced. Moreover, in case $r=1$ the convergence order $\mathcal{O}(n^{-1})$ in (\ref{teq2}) improves the rate  $\mathcal{O}(\ln n \, n^{-1})$ given in \citet{BKZ16}.\\[0.5ex]

\noindent 
{\it Proof of Theorem~\ref{th2}:} We use Theorems \ref{pro1} and  \ref{pro2} with  $H(y)=G_{r,r}(y)$, $h_2(y)= h_{2;r}(y)$
$g_n= r (n-1)+1$ and $b=\min\{r\,,\,2\}$ defined in Proposition \ref{p2} . 

It follows from \eqref{NM2a} with $p= 3/2$ that
\begin{equation}\label{NM2} \mathbb{E}\big(N_n(r)\big)^{- 3/2} \leq C(r) 
\begin{cases} n^{- \min\{r, 3/2\}}, &  r \not= 3/2\\ \ln(n) \, n^{- 3/2}, &   r  = 3/2. 
 \end{cases} 
\end{equation}
The conditions (\ref{add2a}) and (\ref{add2b}) follow from \eqref{corNB} with $p= k/2 < r$, k = 0, 1, 2,  respectively $p=-1, 0$, $r>1$.

Now we estimate the integrals \eqref{gnI1} and \eqref{gnI1b} to obtain a bound in  inequality \eqref{eqpb}. Using   \eqref{fq} for $f_1(z)$ and $f_2(z)$ defined in \eqref{eq33} we find
\begin{equation*}\label{I21} 
J_{1;r;n}(x)= \int\nolimits_{1/g_n}^\infty  \frac{|f_1(x\sqrt{y})|}{\sqrt{y}}dG_{r,r}(y) \leq \frac{c_1^*\,r^r}{\Gamma(r)} \int\nolimits_{1/g_n}^\infty y^{r-3/2}dy = \frac{c_1^*\,r^r\,g_n^{r-1/2}}{\Gamma(r)\,(r-1/2)} 
\end{equation*}
for $0<r<1/2$. If $r=1/2$ then with $x^2/(1+x^2) \leq 1$
$$ J_{1;1/2;n}(x) = \frac{x^2 |p_1|}{8\,\pi p_0^2} \int\nolimits_{1/g_n}^\infty e^{- ((1+x^2)/2)\,y} dy \leq \frac{x^2 |p_1|}{8\,\pi p_0^2\,(1+x^2)/2} \leq \frac{2 |p_1|}{8\,\pi p_0^2}.
$$
Consider  the second term in \eqref{gnI1}. Let $r<1$. Using  now     $c_2^* = \sup_z |f_2(z)|$, then
\begin{equation}\label{I2a}
J^*_{1;r;n}(x) \int\nolimits_{1/g_n}^\infty  \frac{|f_2(x\sqrt{y})|}{y}dG_{r,r}(y) \leq \frac{c_2^*\, r^r}{(1-r)\,\Gamma(r)} \,g_n^{1- r}.
\end{equation} 
If $r=1$ we define the polynomial $P_4(z)$ by  $f_2(z) = z\, P_4(z) \, \varphi(z)$ with $z=x\,\sqrt{y}$ \,and put\, $c_4^{*} = \sup_z \{|P_4(z)|\varphi(z/\sqrt{2})\} < \infty$. \, Then\, $|f_2(z)|  \leq c_4^{*} |z|\,\varphi(z/\sqrt{2})$ and using $|x|\,(1 + x^2/4)^{- 1/2} \leq 2$ we obtain
\begin{equation}\label{I2}
J^*_{1;1;n}(x) \leq \frac{c_4^{*}\,|x|}{\sqrt{2\,\pi}} \int\nolimits^{\infty}_{1/g_n} y^{-1/2} e^{- \,(1 + x^2/4)\,y}\,dy  \leq
\frac{c_4^{*}\,|x| \Gamma(1/2)}{\sqrt{2\,\pi}\,(1 + x^2/4)^{1/2}}  \leq c_4^{*}/\sqrt{2}.
\end{equation}   
and for $0< r \leq 1$ uniform in $x$
\begin{equation}\label{I1}
|I_1(x,n)| \leq g_n^{- 1/2}\, J_{1;r;n}(x)  \,{\rm \bf{1}}_{(0 , 1/2]}(b)+ g_n^{- 1}\, J^*_{1;r;n}(x) \leq C_r g_n^{- r}.
\end{equation}

It remains to estimate $I_2(x,n)$ in (\ref{gnI1b}) for $r > 1$.
Integration by parts for Lebesgue-Stieltjes integrals  and (\ref{add2b}i)  lead to
\begin{eqnarray}\label{eqp5aa}|I_2(x,n)| &\leq& \frac{1}{n}(f_1(x/g_n) + f_2(x/g_n)\,|h_{2;r}(1/g_n)| + I_2^*(x,n)\nonumber\\
&\leq& (c_1^*+c_2^*)c_3^*g_n^{- r} + I_2^*(x,n) 
\end{eqnarray}
with 
\begin{equation}\label{eqp5a} I_2^*(x,n) \,=\,\int_{1/g_n}^\infty  
\Big(\frac{|q_1(x\sqrt{y})|}{2\,n\,\sqrt{g_n} y^{3/2}} + \frac{|q_2(x\sqrt{y})|}{8\,n\,g_n y^2}\Big) \,|h_{2;r}(y)|dy,
\end{equation}
where for $k=1, 2$ functions $f_k(z)$ and $q_k(z)$ are bounded, see \eqref{fq}.

Moreover $g_n y^2 \geq \sqrt{g_n}y^{3\,/2}$ for $y \geq 1/g_n$ and $g_n \leq n\,r$ for $r>1$.

 If $1<r<3/2$ with above defined $c_3^*$ we find 
$$|I_2^*(x,n)| \leq  \frac{(4c_1^{**}+c_2^{**})c_3^*}{8\,n\,\sqrt{g_n}} \int\nolimits_{1/g_n}^\infty  y^{r - \,5/2}  dy = \frac{(4c_1^{**}+c_2^{**})\,r\,c_3^*}{8\,\Gamma(3/2-r)}  \,g_n^{- r}.
$$

If $r > 3/2$ with $c_4^* = \frac{\Di r^{r-1}}{\Di 2 \, \Gamma(r)} \,\sup_y \{(e^{-r\,y/2}\,(|y-1|\,|2 - r| +1)\} < \infty $ we obtain
$$
|I_2^*(x,n)| \leq  \frac{(4c_1^{**}+c_2^{**})c_4^*}{8\,n\,\sqrt{g_n}} \int\nolimits_{1/g_n}^\infty  y^{r - 5/2} e^{-r\,y/2}  dy 
\leq  \frac{(4c_1^{**}+c_2^{**})\,r\,c_4^* \Gamma(r-3/2)}{(r/2)^{r - 3/2} g_n^{3/2}} \, .
$$

For $r = 3/2$ the above estimates of $|I_2^*(x,n)|$ lead to  an exponential integral: 
\begin{eqnarray*}
|I_2^*(x,n)| &\leq & \frac{(4c_1^{**}+c_2^{**})c_4^*}{8\,n\,\sqrt{g_n}} \left(\int\nolimits_{1/g_n}^1  y^{-1} dy + \int\nolimits_1^\infty  e^{-3\,y/2}  dy \right)\\
& \leq & \frac{(4c_1^{**}+c_2^{**})\,3\,c_4^* }{16 } \, \left( \ln n +\ln(3/2)+ \frac{2}{3}\,e^{- 3/2} \right)\,g_n^{- 3/2}.
\end{eqnarray*}

In the latter case $|I_2^*(x,n)| \leq C\, g_n^{- 3/2}$  may be obtained with an analogous procedure as for estimating the above integral  $|I_1(x,n)|$ for $r=1$ in \eqref{I2a}. This proof is omitted because  the rate of convergence in Theorem~\ref{th2}, see \eqref{teq2}, is determined by the negative moment \eqref{NM2}, where the term $\ln n$ cannot be omitted.\\[1ex] 

To obtain (\ref{teq3a}) we calculate integrals in (\ref{gn}), which are similar to that in the proof of Theorem 2           in \citet{CMU20}.  
Using formula 2.3.3.1 in \citet[p.~322]{PBM92} with  $\beta > - r$ and $p=r+x^2/2$:
\begin{equation}\label{prud1}  \int^{\infty}_{0} \frac{y^{\beta}}{e^{ x^2y/2}} dG_{r,r}(y) = 
\frac{r^r}{\Gamma(r)} \int^{\infty}_{0}\frac{y^{\beta +r-1}}{e^{(r+x^2/2)y}} dy = 
\frac{ r^{r} \, \Gamma(\beta + r)}{\Gamma(r)\,(r+x^2/2)^{r + \beta}}
\end{equation}
we compute the first integral  in (\ref{gn}) with $\beta=1/2$ in (\ref{prud1}):
\begin{eqnarray*}
\frac{\partial}{\partial x} \int^{\infty}_{0} \Phi(x\sqrt{y})dG_{r,r}(y)& = & \frac{r^r}{\Gamma(r)\,\sqrt{2\, \pi}}   \int\nolimits^{\infty}_{0} \frac{y^{r-1/2}}{e^{(r+ x^2/2) y}}  dy \\
 &=& \frac{ r^{r} \, \Gamma(r+ 1/2)}{\Gamma(r)\,\sqrt{2\, \pi}\,(r+x^2/2)^{r + 1/2}}=s_{2 r}(x).
\end{eqnarray*}
Hence  
$$\int^{\infty}_{0} \Phi(x\sqrt{y})dG_{r,r}(y) =  S_{2 r}(x)\,.$$
For $r>1/2$ we find with $f_1(x)$ defined in  (\ref{eq33}) and $\beta=1/2$ in (\ref{prud1})
$$
\int\nolimits^{\infty}_{0}\frac{f_1(x\sqrt{y})}{\sqrt{y}} dG_{r,r}(y) =  \frac{p_1\, x |x| \, r^r}{4\, p_0^2 \, \Gamma(r) \,\sqrt{2\, \pi}} 
 \int\nolimits^{\infty}_{0} \frac{y^{r-1/2}}{e^{(r+ x^2/2) y}}dy =  \frac{p_1\, x |x|}{4\, p_0^2}  s_{2 r}(x).
$$
For $r>1$ we obtain with $f_2(x)$ from (\ref{eq33}) and  $\beta = -1/2, 1/2, 3/2$ in (\ref{prud1})
\begin{eqnarray}\label{f2}  %
&&\int\nolimits^{\infty}_{0}\frac{f_2(x\sqrt{y})}{y} dG_{r,r}(y) =  \frac{r^r}{4\, \Gamma(r) \,\sqrt{2\, \pi}} 
 \int\nolimits^{\infty}_{0} \left\{
 3\,x\, y^{r-3/2}\right.\nonumber\\
&&\qquad\qquad\qquad\quad \left. +\Big(1 +\frac{p_2 }{6 p_0^3}\Big)x^3 y^{r -1/2}
   - \,\frac{p_1^2}{8 p_0^4} x^5 \,  y^{r + 1/2}\right\}e^{- (r+ x^2/2) y}dy \nonumber\\
&& \qquad = \frac{x \, r^r}{4\, \Gamma(r) \,\sqrt{2\, \pi}}\left\{
 \frac{3\, \Gamma(r-1/2)}{(r+ x^2/2)^{r-1/2}} +  \Big(1 +\frac{p_2 }{6 p_0^3}\Big)\frac{x^2 \Gamma(r+1/2)}{(r+ x^2/2)^{r +1/2}}\right.\nonumber\\
&& \qquad \qquad \qquad \quad \left. 
- \,\frac{p_1^2 x^4 \Gamma(r+3/2)}{8 p_0^4\, (r+ x^2/2)^{r + 3/2}} \right\},\nonumber\\
&& \qquad = \frac{x}{4}\left\{
 \frac{3\, (2 r+ x^2)}{2 r - 1} +  \Big(1 +\frac{p_2 }{6 p_0^3}\Big)\,x^2 
- \,\frac{p_1^2 x^4 (2 r + 1)}{8 p_0^4\, (2 r+ x^2)} \right\}\,s_{2 r}(x).  
\end{eqnarray}

 The integral $\int_{0}^\infty \Phi(x\sqrt{y})dh_{2;r}(y)$ in (\ref{gn})
is  the same as the integral $J_4(x)$ in the proof of Theorem 2 in \citet{CMU20} where is shown:
$$\frac{1}{n} \sup_x\left|\int_{0}^\infty \Phi(x\sqrt{y})dh_{2;r}(y) -  \frac{ (2 - r) x (x^2 +1)}{4 r (2 r - 1)} \,s_{2 r}(x)\right| \leq c(r)\, n^{- r}\,. $$
With  (\ref{f2}) and $|(r n)^{-1} - g_n^{-1}| \leq (1/r) n^{- 2}$  the term by $1/g_n$ in (\ref{teq3a}) follows.
\hfill $\Box$

\section{Laplace Distribution as Limit for Random Sample Median $M_{N_n}$}

Let $Y(s)\in \mathbb{N}$ be discrete Pareto II distributed with parameter $s>0$, having
probability mass and distribution functions 
\begin{equation}\label{eqP1} 
\mathbb{P}(Y(s)=k)=\frac{s}{s+k-1} - \frac{s}{s+k} \quad \mbox{and} \quad \mathbb{P}\big(Y(s)\leq k\big) =\frac{k}{s+k}, \,k \in \mathbb{N},
\end{equation}
which is a particular class of a  general model of discrete Pareto distributions, obtained by discretization continuous Pareto II (Lomax) 
distributions on integers, see \citet{BK14}.     

Now, let $Y_1(s), Y_2(s), . . . $,  be independent random variables with the same distribution (\ref{eqP1}). Define for $n \in \mathbb{N}$ and $s>0$ the random variable
\begin{equation}\label{eqP3} 
N_n(s)= \max_{1\leq j \leq n} Y_j(s) \quad \mbox{with} \quad \mathbb{P}(N_n(s)\leq k)=\left(\frac{k}{s+k}\right)^{n},\quad n \in \mathbb{N}.
\end{equation}
The distribution of $N_n(s)$ is extremely spread out on the positive integers.

In \citet{CMU20} the following Edgeworth expansion was proved:
\begin{proposition}\label{p3} Let the discrete random variable $N_n(s)$ have distribution function (\ref{eqP3}). 
For $y>0$, fixed $s >0$ and all $n \in \mathbb{N}$ then there exists a real number $C_3(s) > 0$ such that
\begin{equation}\label{eq9aL} 
\sup_{y>0}\left|\mathbb{P}\left(\frac{N_n(s)}{n} \leq y\right) - H_s(y) - \frac{h_{2;s}(y)}{n}\right|\leq \frac{C_3(s)}{n^2},
\end{equation}
\begin{equation}\label{eq9aL2} 
H_s(y) =  e^{-s/y} \,\, \mbox{and} \,\, h_{2;s}(y) = s\, e^{-s/y}\,\big(s-1 + 2 Q_1(n\,y)\big)/\big(2\,y^2\big),\,\, y>0\,,
\end{equation}
where $Q_1(y)$ is defined in (\ref{eq9aa}).
\end{proposition}
{\bf Remarks: 1.}  \citet{Ly10} proved a first order  bound in \eqref{eq9aL} for integer~$s \geq 1$
\begin{equation}\label{ly} 
\left|\mathbb{P}\left(\frac{N_n(s)}{n} \leq x\right) - e^{-s/x}\right| \leq \frac{C_s}{n}, \, C_s=\begin{cases} 
8 e^{-2}/3=0.36..., &\!  s=1, n \geq 2 \\[1ex]
2 e^{-2} = 0.27..., &\! s\geq 2 , n \geq 1         \end{cases}.
\end{equation}
In case $n=1$ and $s=1$ we have $\mathbb{P}\left(N_1(1) \leq x\right) =0$ for $0<x<1$ and
$$ \sup\nolimits_{0<x<1}\left|\mathbb{P}\left(N_1(1) \leq x\right) - e^{-1/x}\right| =\sup\nolimits_{0<x<1}e^{-1/x}= e^{-1}=0.367... \,. $$

{\bf 2.} The continuous function $H_s(y) = e^{-s/y} {\rm \bf{1}}_{(0\,,\,\infty)}(y)$ with parameter $s>0$ is the distribution 
function of the inverse exponential random 
variable  $W(s)=1/V(s)$,  where $V(s)$ is  exponentially distributed  with rate parameter $s>0$.
Both $H_s(y)$  and $\mathbb{P}(N_n(s) \leq y)$ are heavy tailed with shape parameter~1. 

Therefore $\mathbb{E}\big(N_n(s)\big) = \infty$ for all $n \in \mathbb{N}$
and $\mathbb{E}\big(W(s)\big) = \infty$. Moreover:\\[0.3ex]
$\bullet$  First absolute pseudo moment $\nu_1 =\int_{0}^{\infty} x \big|d\big(\mathbb{P}\big(N_n(s) \leq n\,x\big) - e^{-s/x}\big)\big| = \infty$,\\
$\bullet$ Absolute difference moment $\chi_u = \int_{0}^{\infty}x^{u-1} \big|\mathbb{P}\big(N_n(s) \leq n\,x\big) - e^{-s/x}\big|dx < \infty$\\
for  $1\leq u <2$. These statements are proved in \citet[Lemma 2]{CMU20}.
On pseudo moments and some of their generalizations   see e.g. \citet[Chapter 2]{CW92}.\\

Next we estimate the negative moment $\mathbb{E}(N_n(s))^{-p}$, $p>0$, for the random sample size $N_n(s)$:

\begin{corollary}\label{cor52}
Let  $r >0$ and $p>0$. Then for all $n\geq 2$ the following  expansions hold for  negative moments:
\begin{equation}\label{c5}
\mathbb{E}(N_n(s))^{-p} = \begin{cases} {\Di  \frac{\Gamma(p+1)}{s^p\,n^p} + \frac{(s-1)\,p\,\Gamma(p+2)}{2\,s^{p+1} n^{p+1}} + R_{1;n}^*,} &0 < p < 1,\\[1.5ex]
 {\Di \frac{\Gamma(p+1)}{s^p\,n^p} \, + \, R_{2;n}^*,} & 1 \leq p < 2,  \\[1.5ex]
 R_{3;n}^*, & p  \geq 2, 
\end{cases} \end{equation}  
where $|R_{k;n}^*| \leq c_k^*(p) \, n^{- 2}$ for some constants $c_k^*(p) < \infty$, \,\,$k=1, 2, 3$.\\[1ex]
\end{corollary} 

\noindent
{\bf Remarks: 1.} The leading terms in \eqref{c5} and the bound \eqref{eq9aL} lead to the estimate
\begin{equation}\label{NM2aL} \mathbb{E}\big(N_n(s)\big)^{- p} \leq C(p) \, n^{- \min\{p, 2\}},
\end{equation}
where for $0< p \leq 2$ the order of the bound is optimal.\\[0.3ex]
{\bf 2.} In \citet[Corollary 3]{Ben20} the expansion \eqref{c5} for $0 < p < 1$ is given with with an additional term at $n^{- p-1}$.\\[0.5ex] 
                       
\noindent
{\it Proof of Corollary~\ref{cor52}:} As in the beginning of the proof of Corollary 4.2  we obtain 
\begin{equation*}
\mathbb{E}(N_n(s))^{-p} 
= \frac{p}{n^p}  \int\nolimits_{1/n}^\infty \frac{1}{x^{p+1}} \left(H_s(x) + \frac{h_{2;s}(x)}{n}\right)dx  + R_1(n) =I_1 + I_2+ I_3  + R_1(n), 
\end{equation*}
where 
$$ I_1 = \frac{p}{n^p}  \int\nolimits_{0}^\infty \frac{e^{- s/x}}{x^{p+1}} dx  +  R_2(n) = \frac{p}{n^p}  \int\nolimits_{0}^\infty
 y^{p-1}\,e^{- s\,y} dy + R_2(n) = \frac{\Gamma(p+1)}{s^p\,n^p} + R_2(n),
$$
$$  
I_2 =  \frac{p\,s\,(s-1)}{2\,n^{p+1}}  \int\nolimits_{0}^\infty \frac{e^{- s/x}}{x^{p+3}} dx + R_3(n) = \frac{(s-1)\,p\,\Gamma(p+2)}{2\,s^{p+1} n^{p+1}} + R_3(n) , 
$$
considering \eqref{eqA6}
$$
I_3=  \frac{p\,s}{n^{p+1}}  \int\nolimits_{1/n}^\infty \frac{e^{- s/x}\,Q_1(n\,x)}{x^{p+3}} dx =  \frac{p\,s}{n^{p+1}} \sum\nolimits_{k=1}^{\infty} \frac{1}{k\,\pi} \int\nolimits_{1/n}^\infty \frac{e^{- s/x}\,\sin(2\,\pi\, k \,n\, x)}{x^{p+3}} dx
$$
and with \eqref{eq9aL}  of Proposition \ref{p3} 
\begin{equation*}
|R_1(n)| \leq \frac{p}{n^p}  \int\nolimits_{1/n}^\infty \frac{1}{x^{p+1}} \left|\mathbb{P}\left(\frac{N_n(s)}{n} < x\right) - H_{s}(x) - \frac{h_{2sr}(x)}{n}\right|dx \leq \frac{C_3(s)}{n^2}.
\end{equation*}
Since for $\alpha >0$\,  $0 < \beta\leq 2$ and $s>0$
\begin{equation}\label{corL}
 \alpha\int\nolimits_{0}^{1/n} \frac{e^{- s/x}}{x^{\alpha+1}} dx \leq e^{-s\,n} n^\alpha  \leq C(\alpha, \beta, s)\,n^{- \beta}, \quad C(\alpha, \beta, s) = \Big(\frac{\alpha + \beta}{s\,e}\Big)^{\alpha + \beta}
\end{equation} 
we find with $c_2 =p\,C(p,2-p,s)$ \, and \, $c_3 = s p |s-1|\,C(p,2-p,s)/2$
\begin{equation*}
R_2(n) =\frac{p}{n^p}  \int\nolimits_{0}^{1/n} \frac{e^{- s/x}}{x^{p+1}} dx \leq \frac{c_2}{n^2} \,\,\, \mbox{and}\,\,\,
|R_3(n)| = \frac{p\,s\,|s-1|}{2\,n^{p+1}}  \int\nolimits_{0}^{1/n} \frac{e^{- s/x}}{x^{p+3}} dx \leq \frac{c_3}{n^2}.
\end{equation*}
Both remainders decrease 
exponentially with order $n\,e^{- s\,n}$ respectively $n^3\,e^{- s\,n}$.\\
It remains to estimate $I_3$.  Partial integration leads to
$$ 
\int\nolimits_{1/n}^\infty \frac{\sin(2\,\pi\, k \,n\, x)}{e^{s/x}\,x^{p+3}} dx=\frac{n^{p+3} }{e^{s\,n} 2\,\pi\, k \,n} - \int\nolimits_{1/n}^\infty \left(\frac{p+3}{x^{p+4}}+\frac{s}{x^{p+5}}\right)\frac{\,\cos(2\,\pi\, k \,n\, x)}{e^{s/x} 2\,\pi\, k \,n} dx.
$$
Considering \eqref{corL} and $\sum_{k=1}^\infty k^{- 2} = \pi^2/6$ we obtain $|I_3| \leq c(s,p) n^{- 2}$.\hfill $\Box$\\[1ex]

For an asymptotically normally distributed statistic $T_m$ the limit distribution of the standardized $T_{N_{n}(s)}$
is Laplace distribution  $L_{1/\sqrt{s}}(x)$   having density (\ref{eq6LA}) with $\mu= 1/\sqrt{s}$, 
therefore $l_{1/\sqrt{s}}(x) = \sqrt{s/2}\, e^{- \sqrt{2\,s}\,|x|}$. See \citet{BK08} or \citet{ST16}.
%
\begin{theorem}\label{th2s}
Let $s >0$.
Consider the statistic $M_{N_n}$ with random sample size $N_n = N_n(s)$ having distribution function (\ref{eqP3}). 
If for the statistic $M_m(X_1,...X_m)$ inequality (\ref{eq21c}) holds   and $ g_n=n$, then there exists
 a constant $C_s$ such that  $\forall n \in \mathbb{N}$ 
\begin{equation}\label{teq2s}
\sup\nolimits_{x}\Big| \mathbb{P}\Big(2 p_0 \sqrt{n N_n^*/N_n}(M_{N_n} - \theta)\leq  x\Big) - L_{1/\sqrt{s}}(x ; n)\Big| \leq  C_s \, n^{- 3/2} \,,
\end{equation} 
 where   $N_n^*$ is defined in (\ref{Nn}) and
\begin{eqnarray}\label{teq3s}
 L_{1/\sqrt{s}}(x ; n)& = & L_{1/\sqrt{s}}(x) +\, l_{1/\sqrt{s}}(x)\left\{A_{1;s}(x)\, n^{- 1/2} + A_{2;s}(x)\, n^{- 1}\,\right\},\\
  A_{1;s}(x)& = &\frac{p_1\,x\ |x|}{4\,p_0^2} \qquad and \nonumber\\
A_{2;s}(x) &=& \frac{(4-s) x (1+\sqrt{2 s} |x|)}{8 \, s}
 +\, \frac{x^3}{4} \Big(1  \,+\, \frac{p_2}{6 p_0^3}\Big) - \frac{  p_1^2 x^3 |x| \sqrt{2 s}}{32 \, p_0^4} .\qquad \quad
\end{eqnarray} 
\end{theorem}

\noindent
{\bf Remark:} Under the condition (\ref{ly})  a first order expansions  was announced in the conference paper \citet[Theorem 4.1]{BKZ16}. \\[1ex]

\noindent
{\it Proof of Theorem~\ref{th2s}:} We use Theorems \ref{pro1} and \ref{pro2} with   $H(y)=H_{s}(y)$ and $h_2(y)= h_{2;s}(y)$ 
defined in (\ref{eq9aL2}), $b=2$ and $g_n = n$. 

Considering \eqref{corL}  the functions $H_s(1/n)$, $h_{2;s}(1/n)$ and the corresponding integrals decrease even exponentially with order $n\,e^{- s\,n}$ or $n^2\,e^{- s\,n}$, $s>0$. Moreover, $h_{2;s}(0) = \lim_{y \downarrow 0} h_{2;s}(y) =0$.
 Hence conditions (\ref{add2a}) and (\ref{add2b}) are fulfilled.

It remains to estimate $I_2(x,n)$ given in (\ref{gnI1b}). Changing only $h_{2;r}(y)$ by $h_{2;s}(y)$ in the estimations (\ref{eqp5aa}) and (\ref{eqp5a})  of the corresponding  $I_2(x,n)$ in the proof of Theorem \ref{th2},  using partial integration,  the relations \eqref{eq9aL2}, \eqref{fq} and $n y^2 \geq \sqrt{n}y^{3/2}$ for $y \geq 1/n$, then we obtain
\begin{equation*} |I_2(x;n)| \leq (c_1^* + c_2^*)\,c_4(s) n^{- 2} + \frac{(4\,c_1^{**}+c_2^{**})\,(|s-1|+1)\, \Gamma(5/2)}{16\, s^{3/2} n^{3/2}} = c_5(s) n^{- 3/2} .
\end{equation*}

To obtain (\ref{teq3s}) we calculate integrals in (\ref{gn}) for $b = 3/2$ as in the proof of Theorem 5 in \citet{CMU20}.
Here we use formula 2.3.16.3 in \citet[p.~344]{PBM92} with 
$p=x^2/2>0$, $s>0$, $m = 0, 1, 2$:
\begin{equation}\label{prudn2}%
\hspace*{-0.2cm}\int\nolimits^{\infty}_{0} \frac{e^{- x^2 y/2}}{\sqrt{2\pi}\,y^{m -3/2}} d H_s(y) =
\int\nolimits^{\infty}_{0} \frac{s \,e^{- x^2 y/2 - s/y}}{\sqrt{2\pi}\,y^{m +1/2}} dy 
= \left(-1\right)^m \frac{s}{|x|}\frac{\partial^m}{\partial s^m}e^{-\sqrt{2\,s}|x|},
\end{equation}
where 
$$(-1) \frac{s}{|x|} \frac{\partial}{\partial s}e^{-\sqrt{2\,s}|x|} = l_{1/\sqrt{s}}(x) \,\,\, \mbox{and} \,\,\,
\frac{s}{|x|} \frac{\partial^2}{\partial s^2}e^{-\sqrt{2\,s}|x|} = \Big(\frac{|x|}{\sqrt{2\,s}} + \frac{1}{2 \,s} \Big)l_{1/\sqrt{s}}(x).
$$
In the mentioned proof we obtained with (\ref{prudn2}) for $m=1$
$$ \int^{\infty}_{0} \Phi(x\sqrt{y}) d H_s(y) = L_{1/\sqrt{s}}(x)$$
and with (\ref{prudn2}) for $m=2$
$$\sup_x\left|\int_{0}^{\infty}  \Phi(x\sqrt{y})  d h_{2;s}(y)- \frac{(1-s) x (1+\sqrt{2 s} |x|)}{8 \, s} l_{1/\sqrt{s}}(x)\right|
\leq c(s) n^{- 1/2}.$$
Moreover, using (\ref{prudn2}) for $m=1$, we find
$$   \int\limits_{0}^{\infty} \frac{f_1(x\sqrt{y})}{\sqrt{y}} dH_s(y) =\frac{p_1 x |x| s}{4 p_0^2 \sqrt{2\,\pi}} 
\int\limits_{0}^{\infty} \frac{e^{- (x^2/2)y - s/y}}{y^{3/2}} dy = \frac{p_1 x |x|}{4 p_0^2} \, l_{1/\sqrt{s}}(x)
$$
and, finally, with (\ref{prudn2}) for $m=0, 1, 2$, we calculate
\begin{eqnarray*}
\int\nolimits_{0}^{\infty} \frac{f_2(x\sqrt{y})}{y} dH_s(y)& = & \frac{s}{4\,\sqrt{2\,\pi} }\int\nolimits_{0}^{\infty}
\left\{3 x \sqrt{y} + \Big(1 + \frac{p_2}{6\,p_0^3}\Big)x^3 \,y^{3/2}\right.\\
&&\quad \left. - \frac{p_1^2}{8\, p_0^4}x^5\,y^{5/2}\right\} \frac{e^{- (x^2/2)y - s/y}}{y^3} dy\\
&& \hspace*{-2.5cm}  = \,\,\frac{x}{4}\left\{3\, \frac{1+ \sqrt{2\,s}}{2\,s} + \Big(1 + \frac{p_2}{6\,p_0^3}\Big)x^2 - 
\frac{p_1^2\,x^2 \, |x|\,\sqrt{2\,s}}{8\, p_0^4}\right\}\, l_{1/\sqrt{s}}(x).
\end{eqnarray*}

Together with $\mathbb{E} \left(N_n(s)\right)^{- 3/2} \leq C(s) n^{- 3/2}$ for all $s \geq s_0 >0$, proved in 
 \citet[Lemma~3]{CMU20} we proved (\ref{teq2s}). \hfill $\Box$

\section{Cornish-Fisher Expansions for Quantiles of $M_{m}$ and $M_{N_n}$} 
In statistical inference it is of fundamental importance to
obtain the quantiles of the distribution of statistics under consideration.
Statistical applications and modeling with quantile functions are discussed extensively by \citet{Gi00}.
There are very few quantile functions which can be expressed in closed form. The Cornish-Fisher expansions
provide  tools to approximate the quantiles of  probability distributions.

Let $F_{n}(x)$ be a distribution function admitting a Chebyshev-Edgeworth expansion 
in powers of $g_n^{-1/2}$ with $0<g_n \uparrow \infty$ as $n \to \infty$:
\begin{equation}\label{f1}
F_n(x)=G(x)+g(x) \left( a_{1}(x)g_n^{-1/2}+ a_2(x)g_n^{-1} \right)+ R(g_n),\,R(g_n) =  \mathcal{O}(g_n^{-3/2}),
\end{equation}%
where $g(x)$ is the density of a three times  differentiable limit distribution $G(x)$.

\begin{proposition}\label{pro6} 
Let $F_n(x)$ be given by (\ref{f1}) and let $x(u)$ and $u$ be  quantiles of distributions $F_{n}$
and $G$ with the same order $\alpha$, i.e. $ F_n(x(u))=G(u)=\alpha$. Then the following relation holds for $n \to \infty$:
\begin{equation}\label{trivial}
 x(u)=u+ b_{1}(u) g_n^{- 1/2} + b_{2}(u) g_n^{-1} + R^*(g_n),\,\,\,R^*(g_n) =  \mathcal{O}(g_n^{- 3/2}),
\end{equation}
with
\begin{equation}\label{trivial2}
 b_{1}(u)=- a_{1}(u) \quad \mbox{and} \quad b_2(u) =\frac{g'(u)}{2\,g(u)} a_1^2(u)  + a_1'(u) a_1(u) - a_2(u)\,.
 \end{equation}
\end{proposition}
Proposition~\ref{pro6} is  a direct consequence of more general statements, see e.g. \citet[p.~311-315]{Ul11},
\citet[Chapter 5.6.1]{FUS10} or \citet{UAF16} and the references therein.

First we consider random median $M_{N_n}$ if  sample size $N_n=N_n(r)$ is negative binomial distributed with probability 
mass function (\ref{eq5})  and  Student's $t$-distribution $S_{2r}(x)$  is the limit law.
The  second order  expansion (\ref{teq2}) in Theorem~\ref{th2} admits a relation like  (\ref{f1}) with $g_n = r (n-1) +1$ and $a_k(x) = A_{k;r}(x)$, $k=1,2$.
The transfer Proposition~\ref{pro6} implies the following statement:

\begin{corollary}\label{cor61} 
Suppose $r>0$. Let $x=x_\alpha$ and $u=u_\alpha$ be $\alpha$-quantiles of standardized statistic $\mathbb{P}\Big(2 p_0 \sqrt{g_n}(M_{N_n(r)} - \theta) \leq x\Big)$ 
and of the limit Student's $t$-distribution $ S_{2r}(u)$,  respectively. 
Then  with previous definitions the following Cornish-Fischer expansion holds as  $n \to \infty$:
\begin{gather*}\label{eqCF1} 
x=u - \frac{p_1 u |u|}{4 p_0^2 \sqrt{g_n}}{\rm \bf{1}}_{(1/2,\infty)}(r) + \frac{B_2(u)}{g_n}  {\rm \bf{1}}_{(1,\infty)}(r) +
 \left\{\begin{array}{lr}\mathcal{O}(n^{- \min\{r,3/2\}}), &\!\! r\not= 3/2\\ \mathcal{O}(\ln(n)  n^{- 3/2}), &\!\!  r  = 3/2 
 \end{array} \right.,
\end{gather*}
where \quad 
$B_2(u) =\,  \frac{\Di p_1^2 \, u^3}{\Di 8 \, 
 p_0^4} - \frac{\Di (5 - r) \,u^3 \, +\,(5 r +2) u)}{\Di 4(2\,r\, -\,1)} -
\frac{\Di u^3}{\Di 4} \Big(1 + \frac{\Di p_2}{\Di 6 \, p_0^3}\Big).$
\end{corollary}

Next we study the approximation of quantiles for the random mean $M_{N_n}$ if  sample size $N_n=N_n(s)$ is 
based on discrete Pareto distributions with probability 
mass function (\ref{eqP3})  and  Laplace distribution $L_{1/\sqrt{s}}(x)$   is the limit law. Relation
(\ref{teq2s}) in Theorem~\ref{th2s} admits a expansion like  (\ref{f1}) with $g_n = n$ and $a_k(x) = A_{k;s}(x)$, $k=1,2$.
 The transfer Proposition~\ref{pro6} leads now to:
\begin{corollary}\label{cor62} 
Suppose $s >0$. Let $x=x_\alpha$ and $u=u_\alpha$ be $\alpha$-quantiles of standardized statistic 
$\mathbb{P}\Big(2 p_0 \sqrt{n}(M_{N_n(s)} - \theta)\leq x\Big)$ 
and of the limit Laplace distribution $ L_{1/\sqrt{s}}(u)$,  respectively. 
Then  with previous definitions the following Cornish-Fisher expansion holds
\begin{gather*}\label{eqCF1s} 
x=u - \frac{p_1 u |u|}{4 p_0^2 \sqrt{n}} + \frac{B_2(u)}{n}   +
\mathcal{O}(n^{- 3/2}), \quad   \mbox{as} \quad  n \to \infty,
\end{gather*}
 where \quad
$B_2(u) = \frac{\Di p_1^2 \, u^3}{\Di 8 \, p_0^4} \,+\,  \frac{\Di (4-s)\, u \,(1+\sqrt{2 s} |u|)}{\Di 8 \, s}
- \frac{\Di u^3}{\Di 4}\Big(1 +
 \frac{\Di p_2}{\Di 6 p_0^3}\Big).$
\end{corollary}

For the sake of completeness let us consider the  Cornish-Fischer expansion for the median $M_m$, too. Using 
(\ref{eq21c}) with $a_k(x) = f_k(x)$, $k=1,2$, defined in~(\ref{eq33}).
\begin{corollary}\label{cor63}
Let $x=x_\alpha$ and $u=u_\alpha$ be $\alpha$-quantiles of the standardized statistic $\mathbb{P}_{\theta}\Big(2 p_0 \sqrt{2 [m/2]} (M_m -\theta)\leq x\Big)$ 
and of the limit normal distribution $ \Phi(u)$,  respectively. 
Then  with previous definitions the classical Cornish-Fischer expansion holds as  $m \to \infty$:
\begin{gather*}\label{eqCF3} 
x=u - \frac{p_1 u |u|}{4 p_0^2 \sqrt{m}} + \frac{1}{m}\left( \frac{\Di p_1^2 \, u^3}{\Di 8 \,  p_0^4} - \frac{\Di u}{\Di 4} -
\frac{\Di u^3}{\Di 4} \Big(1 + \frac{\Di p_2}{\Di 6 \, p_0^3}\Big)\right)  +
 \mathcal{O}(m^{- 3/2}).
\end{gather*}
\end{corollary}

\section{Acknowledgement}

Proposition \ref{p1}, Theorems \ref{pro1} and \ref{th2} and  Corollary \ref{cor42}  have been obtained under support of the RSF Grant No. 18-11-00132. The paper was prepared within the framework of the Moscow Center for Fundamental and Applied Mathematics, Moscow State University and HSE University Basic Research Programs.  


\end{document}